\documentclass[10pt,reqno]{amsart}
\usepackage{amsmath,amssymb,amsthm,mathrsfs}
\usepackage{graphicx}
\usepackage{hyperref}
\usepackage{indentfirst}
\usepackage{titlesec}
\hypersetup{ colorlinks=false,
  citebordercolor={0 0.7 0},
  linkbordercolor={1 0 0},
  urlbordercolor={0 0 0},
  filebordercolor={0 0 0},
  pdfborder={0 0 0.7},pdftitle={Discrete History Operators: Exact State-Space Realization for Sequence Memory},pdfauthor={Bin Zhou, Rong Gao, Zhao-Yan Li}}
\newtheorem{theorem}{Theorem}[section]
\newtheorem{definition}[theorem]{Definition}
\newtheorem{example}[theorem]{Example}
\newtheorem{lemma}[theorem]{Lemma}
\newtheorem{proposition}[theorem]{Proposition}
\newtheorem{corollary}[theorem]{Corollary}
\newcommand{\R}{\mathbb R}
\newcommand{\C}{\mathbb C}
\newcommand{\N}{\mathbb N}
\newcommand{\F}{\mathbb F}

\newcommand{\cE}{\mathcal E}
\newcommand{\cV}{\mathcal V}
\newcommand{\cF}{\mathcal F}
\newcommand{\cK}{\mathcal K}

\newcommand{\cS}{\mathcal S}
\newcommand{\cC}{\mathcal C}
\newcommand{\cH}{\mathcal H}
\newcommand{\cM}{\mathcal M}
\newcommand{\cN}{\mathcal N}

\newcommand{\rank}{\operatorname{rank}}
\newcommand{\row}{\operatorname{row}}

\newcommand{\T}{\mathrm T}
\newcommand{\Hh}{\mathrm H}

\newcommand{\vPi}{\varPi}

\newcommand{\norm}[1]{\left\lVert #1\right\rVert}
\newcommand{\abs}[1]{\left\lvert #1\right\rvert}

\renewenvironment{abstract}{%
  \par\addvspace{0.9cm}
  \begin{center}\bfseries Abstract\end{center}
  \vspace{0.3em}\normalfont\normalsize\par\indent
}{\par}
\renewenvironment{keywords}{%
  \par\vspace{1.0em}\noindent\textbf{Keywords.}\enspace
}{\par}

\renewcommand{\footnoterule}{%
  \kern -3pt
  \hrule width 0.45\textwidth
  \kern 3pt
}

\titleformat{\section}[hang]
  {\normalfont\Large\bfseries\raggedright}
  {\thesection.}
  {0.8em}
  {}

\titlespacing*{\section}
  {0pt}
  {2ex plus 1ex minus .2ex}
  {1.2ex}

\makeatletter
\newcommand{\fundingnote}[1]{%
  \begingroup
  \def\@makefnmark{}%
  \def\@thefnmark{}%
  \@footnotetext{%
    \normalfont\footnotesize
    #1%
  }%
  \endgroup
}
\makeatother

\makeatletter
\renewcommand{\subsection}{%
  \@startsection{subsection}{2}%
  \z@
  {.5\linespacing\@plus.7\linespacing}%
  {.5\linespacing}%
  {\normalfont\bfseries}%
}
\makeatother

\begin{document}
\begin{center}
{\LARGE\bfseries
Discrete History Operators:\\[3pt]
Exact State-Space Realization for Sequence Memory\par}
\vspace{1.2cm}
{\large Bin Zhou$^{a,*}$, Rong Gao$^b$, Zhao-Yan Li$^{b,c}$\par}
\vspace{0.65cm}
{\normalsize
$^a$ Center for Control Theory and Guidance Technology,\\
Harbin Institute of Technology, Harbin 150001, China\\
National Key Laboratory of Complex System Control and Intelligent Agent Cooperation\\[4pt]
$^b$ School of Mathematics, Harbin Institute of Technology,\\
Harbin 150001, China\\[4pt]
$^c$ Zhengzhou Research Institute of HIT, Zhengzhou 450000, China\par}
\end{center}

\fundingnote{%
This paper was supported in part by
the National Natural Science Foundation
of China for Distinguished Young Scholars
under Grant 62125303;
by the National Natural Science Foundation
of China--``Qisun Ye'' Science Foundation
under Grant U2441243;
by the National Natural Science Foundation
of China under Grants 62573165 and 62521005;
the Natural Science Foundation of
Heilongjiang Province under Grant
FG2025F001;
and by the Fundamental and Interdisciplinary
Disciplines Breakthrough Plan of the
Ministry of Education of China
JYB2025XDXM206.
$^*$Corresponding author:
Bin Zhou.
Email: \texttt{binzhou@hit.edu.cn}.
}

\begin{abstract}
A prescribed finite-dimensional analysis of an
 evolving discrete history need not admit an
 exact recurrence driven only by its current
 state and the entering input. We study exact
 realizations of causal linear discrete history
 operators (DHOs). Operator factorization
 and kernel inclusion characterize exact
 closure, while the cross-projection defect
 characterizes closure failure for
 orthonormal analyses.
For finite-dimensional history spaces
 over fixed finite horizons, we distinguish
 two minimality problems.
Reachable task quotients determine the
 minimum state dimension for preserving
 prescribed future tasks on reachable
 histories. Backward invariant envelopes
 characterize the minimum dimension of
 exact linear realizations recovering
 prescribed analyses on full history
 spaces within the admissible class.
The minimum number of additional
 independent state directions is given
 by the rank of the propagated closure
 defect operator restricted to the
 kernel of the current analysis.
For geometrically weighted countable
 histories, matrix mode transport yields
 structured exact recurrences and a
 uniform exponential decay estimate
 for state transitions.
Additional conditions provide uniform
 Gram bounds and refined decay estimates.
A stationary Meixner specialization yields
 explicit recurrence matrices, an $O(N)$
 arithmetic state update, and
 minimum augmentation dimensions over
 finite horizons for selected coefficients.
Numerical experiments examine
 recurrence accuracy, state transition
 bounds, Gram conditions, and
 minimum realization dimensions.
\end{abstract}

\section{Introduction}
In estimation and monitoring of
discrete-time systems, a required
quantity may depend on the history
of sampled inputs or measurements.
Examples include filtered signals
and moments over sliding windows
\cite{chen2015sliding,kailath1980linear}.
An exact recursive implementation
must reproduce the specified quantity
as each new sample arrives.

In this paper, we study prescribed
linear analyses of discrete histories.
Let $h_k\in\cH_k$ denote the available
sample history at time $k$, and let
$x_k=\cK_kh_k\in\F^{N_k}$
be its prescribed linear analysis.
We refer to $\cK_k$ as a discrete
history operator (DHO).
Suppose the history evolves according
to
$h_{k+1}=\cS_kh_k+\cE_{k+1}u_{k+1}$,
where $\cS_k$ transports the existing
history and $\cE_{k+1}$ inserts the
new input sample.
From a systems viewpoint, this is
a linear time-varying history system
with prescribed output $x_k$.
The realization problem is to determine
whether this output contains enough
information for its own exact
recursive update. If not, we seek
the minimum additional state information
needed to preserve the prescribed
output under future history transport.
Without an exact recurrence,
direct evaluation of the prescribed
output $x_k=\cK_kh_k$ may require
retaining and revisiting the full history.

An exact recurrence, when it exists,
has the form
\begin{equation}
x_{k+1}=A_kx_k+B_{k+1}u_{k+1}.
\label{eq:ssm}
\end{equation}
The coefficient matrices are required
to be independent of the realized
inputs, so \eqref{eq:ssm} is a discrete
linear time-varying (DLTV)
state-space model.
Exact closure holds precisely when
there exist matrices $A_k$ satisfying
$\cK_{k+1}\cS_k=A_k\cK_k$,
with $B_{k+1}=\cK_{k+1}\cE_{k+1}$.
In this case, the prescribed output
can itself serve as a recursively
updated state. The general closure
conditions are developed in
Section~\ref{sec:general discrete}. When exact closure fails, we seek
an augmented state
$\hat x_k=\hat{\cK}_kh_k$
that admits an exact recurrence
and recovers the prescribed output
through $x_k=F_k\hat x_k$.
For this realization problem,
the recovery and recurrence
identities must hold on the full
history spaces, not merely along
reachable trajectories.
Thus the prescribed output is fixed,
whereas the realization state may
be redesigned.

With $x_0=0$, repeated substitution
in \eqref{eq:ssm} gives
\begin{equation}
x_k=\sum_{j=1}^{k}\cF(k,j)B_ju_j,
\qquad
\cF(k,s):=
\begin{cases}
A_{k-1}\cdots A_s, & k>s,\\
I, & k=s,
\end{cases}
\quad k\ge s\ge0.
\label{eq:ssm-kernel}
\end{equation}
Thus, under zero initialization,
every such recurrence induces
a linear analysis of the input
history, with kernel
$\cF(k,j)B_j$.
The zero initial state is used only
for this representation.
A given initial history $h_s$
may instead determine the state
through $x_s=\cK_sh_s$.
Conversely, a prescribed history
analysis need not admit an exact
closed recurrence without
augmentation.

A simple filtered output illustrates
this obstruction.
Consider $x_k=u_k+u_{k-1}$ for
$k\ge1$.
The histories
$[u_{k-1},u_k]=[0,1]$ and $[1,0]$
both give $x_k=1$.
For the same entering input
$u_{k+1}=0$, however, the next
prescribed outputs are
\[
x_{k+1}
=
\begin{cases}
1, & [u_{k-1},u_k]=[0,1],\\
0, & [u_{k-1},u_k]=[1,0].
\end{cases}
\]
Hence the prescribed output alone
cannot determine its next value
for every input history.
Introducing
$\hat x_k=[x_k,u_k]^{\T}$
gives the exact realization
\[
\hat x_{k+1}
=
\begin{bmatrix}
0 & 1\\
0 & 0
\end{bmatrix}\hat x_k
+
\begin{bmatrix}
1\\
1
\end{bmatrix}u_{k+1},
\]
with $x_k=[1,0]\hat x_k$.
Within the admissible class of
linear realizations, one additional
state coordinate is both necessary
and sufficient.
This example illustrates minimum
exact augmentation for preserving
a prescribed output functional.

Classical realization theory studies
state-space representations of
specified input--output behavior,
including controllability, observability,
and minimality
\cite{ho1966realization,
kailath1980linear,Kalman1963}.
Related geometric methods use
factorization and invariant subspaces
\cite{arbib1975timevarying,
weiss1972controllability,
wonham1985geometric}.
Behavioral system theory provides
another framework for linear system
representations
\cite{willems1986partI,
willems1986partII}.
Projection-based model reduction
seeks to reproduce or approximate
selected aspects of system behavior
\cite{antoulas2005approximation,
benner2015survey}.

For linear time-invariant systems,
minimum realizations of input--output
behavior under zero initial conditions
are controllable and observable.
By contrast, exact recovery on full
history spaces may require retaining
directions outside the reachable
history space when they affect future
prescribed analyses.
When the future tasks coincide with
the prescribed analyses, the minimum
dimensions on reachable and full
history spaces are determined by the
ranks of the same future task operator
on these spaces, respectively.
The propagated closure defect rank
further identifies the additional
independent directions beyond the
current analysis, accounting for
overlapping information across time.

Recurrent and structured state-space
models provide finite-dimensional
representations of sequence histories
\cite{gu2022s4,gupta2022dss,
orvieto2023lru}.
Continuous projection methods,
including HiPPO and Legendre memory
units, provide related constructions
for continuous histories
\cite{gu2020hippo,gu2023trainhippo,
voelker2019lmu}.
For discrete histories, Lag Operator
SSMs derive updates from evolving
bases and lag operators
\cite{tomonaga2026lag}.
SaFARi develops state-space
constructions using general
frames and bases
\cite{babaei2025safari}.
Related work includes orthogonal
rational function representations
\cite{vangucht2003orf},
recursive discrete moments
\cite{chen2015sliding,
wang2006tchebichef},
and updating and downdating
of orthogonal polynomials
\cite{elhay1991update}.
Meixner-like filters have been used
in recursive realizations
\cite{maraoui2017meixner},
while Meixner polynomials have been
applied to model reduction
\cite{xu2024meixner}.
These constructions do not by
themselves provide a general
exact closure criterion or
a minimum augmentation rule for
an independently prescribed
discrete history analysis.

The main contributions are as follows.
\begin{enumerate}
\item We apply operator factorization
and kernel inclusion to characterize
exact closure of prescribed discrete
history analyses. For orthonormal
analyses, we derive a cross-projection
defect identity describing the failure
of projected history transport
to reproduce the prescribed analysis
of the updated history.

\item We distinguish two finite-horizon
minimality problems.
Reachable task quotients determine
the minimum state dimension for
preserving prescribed future tasks
on reachable histories.
For exact recovery of prescribed
analyses on full history spaces,
backward invariant envelopes determine
the minimum state dimension within
the admissible class, which allows
state coordinates to be redesigned.
The minimum number of additional
independent state directions equals
the rank of the propagated closure
defect operator restricted to the
kernel of the current analysis.

\item We construct exact realizations
for geometrically weighted countable
histories using matrix mode transport
with possibly noncommuting matrices.
We establish a uniform exponential
decay estimate for state transitions
and sufficient conditions for uniform
Gram bounds and refined decay estimates.
A stationary Meixner specialization
provides explicit recurrence matrices,
an $O(N)$ arithmetic state update,
and exact minimum augmentation
dimensions over finite horizons
for selected coefficients.
\end{enumerate}

Section~\ref{sec:notation} introduces
the notation and preliminaries.
Section~\ref{sec:general discrete}
develops exact closure conditions,
cross-projection defects, and
reachable task quotients.
Section~\ref{sec:timevaring}
studies minimum exact augmentation
through backward invariant envelopes
and ranks of propagated closure
defects.
Section~\ref{sec:structured-realizations}
presents structured exact realizations
for countable histories.
Section~\ref{sec:numerical} reports
numerical experiments, and
Section~\ref{sec:conclusions}
concludes the paper.
\section{Notation and preliminaries}
\label{sec:notation}

The symbols $\R$, $\C$, $\mathbb Z$, and
$\N$ denote the real numbers, complex
numbers, integers, and nonnegative
integers, respectively.
The scalar field $\mathbb F$ is either
$\R$ or $\C$, with $\mathbb F^0=\{0\}$.
The superscripts $\T$ and $\Hh$ denote
transpose and conjugate transpose, while
$^*$ denotes the adjoint of a Hilbert-space
operator. Inner products on complex
Hilbert spaces are linear in the first
argument and conjugate linear in the second.
For Hilbert spaces $X$ and $Y$,
$\mathcal L(X,Y)$ denotes the space of
bounded linear operators from $X$ to $Y$.
For a linear operator $T$, $\ker T$,
$\operatorname{range}T$, and
$\operatorname{rank}T$ denote its kernel,
range, and rank, respectively.
For a matrix $M$, $\operatorname{row}M$
denotes its row space, and
$\sigma_{\min}(M)$ its smallest singular
value. For a subspace $W$ of a
finite-dimensional vector space $V$,
$\operatorname{codim}_V W=\dim V-\dim W$.
The notations $\operatorname{col}(\cdot)$,
$\operatorname{diag}(\cdot)$, and
$\operatorname{tr}(\cdot)$ denote vertical
concatenation, diagonal matrix formation,
and the matrix trace, respectively.
The symbol $I_n$ denotes the $n\times n$
identity matrix, while $I$ denotes
an identity matrix or operator whose
dimension or underlying space is clear
from context.
For a Hermitian matrix $H$,
$\lambda_{\min}(H)$ and $\lambda_{\max}(H)$
denote its smallest and largest eigenvalues.
For Hermitian matrices $X$ and $Y$ of
the same size, $X\succeq Y$ and $X\succ Y$
mean that $X-Y$ is positive semidefinite
and positive definite, respectively.
Unless otherwise stated, $\norm{\cdot}$
denotes the Euclidean norm for vectors,
the spectral norm for matrices, and the
operator norm for bounded Hilbert-space
operators. The notation
$\norm{\cdot}_{\mathrm F}$ denotes the
Frobenius norm for matrices and the
Hilbert--Schmidt norm for operators
between finite-dimensional Hilbert spaces.
For an invertible matrix $M$, its
spectral condition number is
$\kappa_2(M)=\norm{M}\norm{M^{-1}}$.
We use $O(\cdot)$ in the standard
asymptotic sense.

The following elementary rank identity
will be used in the analysis of
minimum exact augmentation.

\begin{lemma}
\label{lem:rank-increment}
Let $V$ be a finite-dimensional vector
space over $\mathbb F$, and let
$K:V\to\mathbb F^p$ and
$L:V\to\mathbb F^q$ be linear maps.
Then
\begin{equation}
\operatorname{rank}
\begin{bmatrix}
K\\L
\end{bmatrix}
-\operatorname{rank}K
=
\operatorname{rank}
\left(L|_{\ker K}\right).
\label{eq:rank-increment}
\end{equation}
\end{lemma}

\begin{proof}
Let $S:V\to\mathbb F^{p+q}$ be defined by
$S(v)=(Kv,Lv)$.
Then $\ker S=\ker K\cap\ker L$.
Applying the rank--nullity theorem
to $S$ and $K$ gives
\[
\operatorname{rank}S-\operatorname{rank}K
=
\dim\ker K-\dim(\ker K\cap\ker L).
\]
Since $\ker(L|_{\ker K})=\ker K\cap\ker L$,
a further application of the same theorem
to $L|_{\ker K}$ yields
$\operatorname{rank}S-\operatorname{rank}K
=\operatorname{rank}(L|_{\ker K})$,
as required.
\end{proof}

\section{General theory of discrete history operators}\label{sec:general discrete}
We characterize exact closure of
prescribed outputs of time-varying
history systems through operator
factorization and kernel inclusion.
For orthonormal analyses, we characterize
closure failure using the
cross-projection defect.
We then characterize minimum state
dimensions for preserving specified
future tasks on reachable histories
over finite horizons.
\subsection{Exact closure and cross-projection defects}
\label{subsec:countable}

Let $\cH_k$ be a Hilbert space over $\F$
representing the history available at time $k$,
and let $h_k\in\cH_k$ denote the corresponding
history.

\begin{definition}\label{def:dho}
A DHO is a family of bounded linear
operators
$\{\cK_k\}_{k\ge0}$ with
$\cK_k\in\mathcal L(\cH_k,\F^{N_k})$,
where $N_k\ge0$ is the number of
prescribed analysis coordinates at
time $k$. For $h_k\in\cH_k$, the
corresponding DHO state is
\begin{equation}
x_k=\cK_kh_k.
\label{eq:dho-state}
\end{equation}
Since $\cH_k$ contains only the history
available at time $k$, the DHO is causal
by construction.
\end{definition}

If
$\cK_kh=
\operatorname{col}(\langle h,\phi_{k,n}\rangle)
_{n=0}^{N_k-1}$
for an orthonormal family
$\{\phi_{k,n}\}_{n=0}^{N_k-1}\subset\cH_k$,
then $\cK_k$ is bounded and surjective, with
$\cK_k^*z=\sum_{n=0}^{N_k-1}z_n\phi_{k,n}$
for $z\in\F^{N_k}$.
Orthonormality gives
$\cK_k\cK_k^*=I_{N_k}$.
Consequently, $\cK_k^*\cK_k$ is the
orthogonal projector onto
$\operatorname{span}\{\phi_{k,n}:0\le n<N_k\}$.
Thus $\cK_kh_k$ gives the projection
coordinates of the available history.

Let $m\ge1$ be the input dimension.
For each $k\ge0$, let
$\cS_k\in\mathcal L(\cH_k,\cH_{k+1})$
be the history transport operator and
$\cE_{k+1}\in
\mathcal L(\F^m,\cH_{k+1})$
be the insertion operator for the
new sample $u_{k+1}\in\F^m$.
The history evolves according to
\begin{equation}
h_{k+1}
=
\cS_kh_k+\cE_{k+1}u_{k+1}.
\label{eq:history-migration}
\end{equation}
All adjoints are taken with respect to
the corresponding Hilbert-space inner
products, with standard inner products
on finite-dimensional input and state
spaces.

We first characterize exact realization
of a prescribed DHO.
Assume that each
$\cK_k\in\mathcal L(\cH_k,\F^{N_k})$
is surjective.

\begin{theorem}
\label{thm:countable-operator-closure}
Matrices
$A_k\in\F^{N_{k+1}\times N_k}$ and
$B_{k+1}\in\F^{N_{k+1}\times m}$
satisfy
\begin{equation}
x_{k+1}
=
A_kx_k+B_{k+1}u_{k+1}
\label{eq:countable-state-recurrence}
\end{equation}
for every $h_k\in\cH_k$ and
$u_{k+1}\in\F^m$ if and only if
\begin{equation}
\cK_{k+1}\cS_k=A_k\cK_k,
\qquad
B_{k+1}=\cK_{k+1}\cE_{k+1}.
\label{eq:operator-closure}
\end{equation}
When such matrices exist, both $A_k$
and $B_{k+1}$ are uniquely determined.
\end{theorem}

\begin{proof}
If \eqref{eq:operator-closure} holds,
applying $\cK_{k+1}$ to
\eqref{eq:history-migration} gives
\[
x_{k+1}
=\cK_{k+1}h_{k+1}
=\cK_{k+1}(\cS_kh_k+\cE_{k+1}u_{k+1})=A_kx_k+B_{k+1}u_{k+1}.
\]
This proves
\eqref{eq:countable-state-recurrence}.

Conversely, suppose
\eqref{eq:countable-state-recurrence}
holds for every $h_k\in\cH_k$ and
$u_{k+1}\in\F^m$.
By \eqref{eq:dho-state} and
\eqref{eq:history-migration},
\[
\cK_{k+1}\cS_kh_k+
\cK_{k+1}\cE_{k+1}u_{k+1}
=
A_k\cK_kh_k+B_{k+1}u_{k+1}.
\]
Setting $u_{k+1}=0$ yields
$\cK_{k+1}\cS_k=A_k\cK_k$,
while setting $h_k=0$ gives
$B_{k+1}=\cK_{k+1}\cE_{k+1}$.
Thus \eqref{eq:operator-closure} holds.

To prove uniqueness, suppose that
$A_k$ and $\hat A_k$ both satisfy
the first identity in
\eqref{eq:operator-closure}.
Then
$(A_k-\hat A_k)\cK_k=0$.
Since $\cK_k$ is surjective, every
$z\in\F^{N_k}$ can be written as
$z=\cK_kh$ for some $h\in\cH_k$.
Consequently,
$(A_k-\hat A_k)z=0$ for every
$z\in\F^{N_k}$, so
$A_k=\hat A_k$.
The second identity in
\eqref{eq:operator-closure} uniquely
determines $B_{k+1}$.
\end{proof}

The factorization condition
\eqref{eq:operator-closure} admits an
equivalent characterization in terms
of the history directions discarded
by the current DHO state.

\begin{theorem}
\label{thm:quotient-parameterization}
Assume that $\cK_k$ is surjective for every
$k\ge0$. The DHO admits the exact recurrence
\eqref{eq:countable-state-recurrence} at every
time step if and only if
\begin{equation}
\ker\cK_k
\subseteq
\ker(\cK_{k+1}\cS_k),
\qquad k\ge0.
\label{eq:kernel-invariance}
\end{equation}
When \eqref{eq:kernel-invariance} holds, let
$q_k:\cH_k\to\cH_k/\ker\cK_k$
be the quotient map and define
\[
U_k(q_kh)=\cK_kh,
\qquad
\widetilde{\cS}_k(q_kh)=q_{k+1}\cS_kh.
\]
Then $U_k$ is an isomorphism from
$\cH_k/\ker\cK_k$ onto $\F^{N_k}$,
$\widetilde{\cS}_k$ is well defined and bounded,
and
\begin{equation}
\cK_k=U_kq_k,
\qquad
A_k=U_{k+1}\widetilde{\cS}_kU_k^{-1}.
\label{eq:quotient-realization}
\end{equation}

Conversely, let $\{\cN_k\}$ be a family of
closed subspaces of $\cH_k$ with
$\dim(\cH_k/\cN_k)=N_k$.
Suppose that
$\cS_k\cN_k\subseteq\cN_{k+1}$.
For any isomorphisms
$U_k^{\cN}:\cH_k/\cN_k\to\F^{N_k}$,
the operators
$\cK_k^{\cN}=U_k^{\cN}q_k^{\cN}$,
where $q_k^{\cN}$ denotes the quotient map,
define a DHO satisfying exact closure.
\end{theorem}

\begin{proof}
Suppose first that the exact recurrence holds.
By Theorem~\ref{thm:countable-operator-closure},
$\cK_{k+1}\cS_k=A_k\cK_k$.
Hence, for every $h\in\ker\cK_k$,
$\cK_{k+1}\cS_kh=A_k\cK_kh=0$,
which proves \eqref{eq:kernel-invariance}.

Conversely, suppose
\eqref{eq:kernel-invariance} holds.
Since $\cK_k$ is bounded,
$\ker\cK_k$ is closed, and
$\cH_k/\ker\cK_k$ is a Hilbert space.
Define $U_k(q_kh)=\cK_kh$.
If $q_kh=q_k\hat h$, then
$h-\hat h\in\ker\cK_k$, so
$\cK_kh=\cK_k\hat h$.
Thus $U_k$ is well defined.
Since $\cK_k(h+v)=\cK_kh$ for every
$v\in\ker\cK_k$, the quotient norm gives
\[
\norm{U_k(q_kh)}
=\norm{\cK_kh}
\le\norm{\cK_k}
\norm{q_kh}_{\cH_k/\ker\cK_k}.
\]
Hence $U_k$ is bounded.
It is injective because
$U_k(q_kh)=0$ implies $q_kh=0$,
and surjective because $\cK_k$ is
surjective. Thus $U_k$ is a bounded
bijection onto $\F^{N_k}$.
Since its domain is finite dimensional,
$U_k^{-1}$ is also bounded.
Therefore $U_k$ is an isomorphism and
$\cK_k=U_kq_k$.

Next define
$\widetilde{\cS}_k(q_kh)=q_{k+1}\cS_kh$.
If $q_kh=q_k\hat h$, then
$h-\hat h\in\ker\cK_k$.
By \eqref{eq:kernel-invariance},
$\cS_k(h-\hat h)\in\ker\cK_{k+1}$,
and hence
$q_{k+1}\cS_kh=q_{k+1}\cS_k\hat h$.
Thus $\widetilde{\cS}_k$ is well defined
and linear. Since both quotient spaces
are finite dimensional, it is bounded.

Define
$A_k=U_{k+1}\widetilde{\cS}_kU_k^{-1}$.
For every $h\in\cH_k$,
\[
A_k\cK_kh
=U_{k+1}\widetilde{\cS}_kq_kh=U_{k+1}q_{k+1}\cS_kh=\cK_{k+1}\cS_kh.
\]
Consequently,
$\cK_{k+1}\cS_k=A_k\cK_k$.
With $B_{k+1}=\cK_{k+1}\cE_{k+1}$,
Theorem~\ref{thm:countable-operator-closure}
gives the exact recurrence and
\eqref{eq:quotient-realization}.

Finally, let $\{\cN_k\}$ satisfy the
stated conditions and define
$\cK_k^{\cN}=U_k^{\cN}q_k^{\cN}$.
These operators are bounded and
surjective, with
$\ker\cK_k^{\cN}=\cN_k$.
The invariance condition implies
\[
\cS_k(\ker\cK_k^{\cN})
=
\cS_k\cN_k
\subseteq
\cN_{k+1}
=
\ker\cK_{k+1}^{\cN}.
\]
Thus \eqref{eq:kernel-invariance} holds
for the constructed family, and the
equivalence established above proves
exact closure.
\end{proof}

Condition~\eqref{eq:kernel-invariance}
requires every direction in
$\ker\cK_k$ to remain invisible
to the next analysis after transport.
Its failure is characterized by the
cross-projection defect for orthonormal
analyses and motivates the minimum
exact augmentation problem.

\begin{proposition}
\label{prop:cross-projection-exactness}
Assume $\cK_k\cK_k^*=I_{N_k}$ and define
$\vPi_k=\cK_k^*\cK_k$.
Let $I_{\cH_k}$ denote the identity operator
on $\cH_k$, and set
\begin{equation}
A_k^{\mathrm{cp}}
=
\cK_{k+1}\cS_k\cK_k^*,
\qquad
\Delta_k^{\mathrm{cp}}
=
\cK_{k+1}\cS_k(I_{\cH_k}-\vPi_k).
\label{eq:cross-projection-defect}
\end{equation}
We call $\Delta_k^{\mathrm{cp}}$
the cross-projection defect.
Then
\begin{equation}
\cK_{k+1}\cS_k
=
A_k^{\mathrm{cp}}\cK_k
+\Delta_k^{\mathrm{cp}}.
\label{eq:cross-projection-identity}
\end{equation}
Consequently, with
$B_{k+1}=\cK_{k+1}\cE_{k+1}$,
the matrix $A_k^{\mathrm{cp}}$ gives
exact closure if and only if
$\Delta_k^{\mathrm{cp}}=0$, equivalently,
\eqref{eq:kernel-invariance}.
Moreover, for every $h\in\cH_k$,
\begin{equation}
A_k^{\mathrm{cp}}\cK_kh
=
\cK_{k+1}\cS_k\vPi_kh.
\label{eq:cross-projection-projected-history}
\end{equation}
For $x_k=\cK_kh_k$, the one-step
state error satisfies
\begin{equation}
x_{k+1}
-
\bigl(A_k^{\mathrm{cp}}x_k
+B_{k+1}u_{k+1}\bigr)
=
\Delta_k^{\mathrm{cp}}h_k.
\label{eq:cross-projection-one-step-error}
\end{equation}
\end{proposition}

\begin{proof}
Since $\cK_k\cK_k^*=I_{N_k}$,
the operator $\vPi_k=\cK_k^*\cK_k$
satisfies
$\vPi_k^2=\vPi_k$ and
$\vPi_k^*=\vPi_k$.
Moreover,
$\operatorname{range}\vPi_k
=\operatorname{range}\cK_k^*$.
Thus $\vPi_k$ is the orthogonal
projector onto
$\operatorname{range}\cK_k^*$, and
\(\operatorname{range}(I_{\cH_k}-\vPi_k)
=
\ker\cK_k.\)
By the definition of $A_k^{\mathrm{cp}}$,
\(A_k^{\mathrm{cp}}\cK_k
=
\cK_{k+1}\cS_k\vPi_k.\)
Decomposing $I_{\cH_k}$ as
$\vPi_k+(I_{\cH_k}-\vPi_k)$
therefore gives
\eqref{eq:cross-projection-identity}.

Since
$\operatorname{range}(I_{\cH_k}-\vPi_k)
=\ker\cK_k$,
the equality $\Delta_k^{\mathrm{cp}}=0$
holds if and only if
$\cK_{k+1}\cS_kh=0$ for every
$h\in\ker\cK_k$.
This is precisely
\eqref{eq:kernel-invariance}, which is
equivalent to exact closure by
Theorem~\ref{thm:quotient-parameterization}.

Applying
$A_k^{\mathrm{cp}}\cK_k
=\cK_{k+1}\cS_k\vPi_k$
to $h\in\cH_k$ proves
\eqref{eq:cross-projection-projected-history}.
Finally, combining
\eqref{eq:cross-projection-identity}
with \eqref{eq:history-migration}
and $B_{k+1}=\cK_{k+1}\cE_{k+1}$
gives
\eqref{eq:cross-projection-one-step-error}.
\end{proof}

Lag Operator SSMs~\cite{tomonaga2026lag}
derive discrete recurrence matrices from inner products
of adjacent bases after replacing the
previous history by its finite-dimensional
approximation.
Their update construction does not by
itself establish exact recovery of an
independently prescribed analysis of
the full discrete history.
In our setting, this distinction is
expressed by the factorization condition
\eqref{eq:operator-closure}, while
\eqref{eq:cross-projection-one-step-error}
identifies the resulting closure error.

SaFARi~\cite{babaei2025safari}
constructs state-space representations
from general frames and bases.
Its truncation analysis distinguishes
representation and mixing errors, and
establishes optimal reconstruction
properties for a fixed frame.
Our minimum augmentation problem instead
requires exact recovery of prescribed
analysis coefficients for every history.
Theorem~\ref{thm:minimal-window-augmentation}
and
Corollary~\ref{cor:propagated-closure-defect-rank}
determine the minimum number of
additional independent state directions
under this constraint.

\subsection{Reachable task quotients}
\label{subsec:reachable-task-quotients}
The preceding results characterize
exact closure of prescribed analyses
on full history spaces.
We now consider minimum state
realizations that preserve specified
future task outputs on reachable
histories over fixed finite horizons.
The task operators are linear maps
of the histories, and the history
spaces considered here are
finite-dimensional.

Fix $s<t$ and assume that $\cH_k$ is
finite dimensional for $s\le k\le t$.
Let $\mathscr R_s\subseteq\cH_s$ be
a prescribed subspace of admissible
initial histories. A history
$h_k\in\cH_k$ is reachable from
$\mathscr R_s$ if it can be generated
by \eqref{eq:history-migration}
from some $h_s\in\mathscr R_s$
under a sequence of entering inputs.
The reachable history space
$\mathscr R_k$ consists of all such
histories and satisfies
\begin{equation}
 \mathscr R_{k+1}
 =
 \cS_k\mathscr R_k
 +
 \operatorname{range}\cE_{k+1},
 \qquad s\le k<t.
 \label{eq:reachable-history-recursion}
\end{equation}

Let $\cC_k\in\mathcal L(\cH_k,\F^{p_k})$
be a prescribed linear task operator,
with $\cC_kh_k$ representing the
quantity of interest extracted from
the history at time $k$.
For $i>k$, define
$\cS(i,k)=\cS_{i-1}\cdots\cS_k$, and set
$\cS(k,k)=I_{\cH_k}$.
Define the finite-horizon future-task operator
$\mathscr O_{t,k}\in
\mathcal L(\cH_k,\F^{p_k+\cdots+p_t})$
by
\begin{equation}
\mathscr O_{t,k}h
=
\operatorname{col}\!\left(
\cC_kh,
\cC_{k+1}\cS(k+1,k)h,
\ldots,
\cC_t\cS(t,k)h
\right).
\label{eq:future-task-operator}
\end{equation}
Thus $\mathscr O_{t,k}h$ collects
the prescribed task values generated
from $h$ over $[k,t]$ under zero input.

\begin{theorem}
\label{thm:finite-horizon-quotient-synthesis}
Define the subspaces $\cM_k$ backward by
\begin{equation}
\cM_t=\mathscr R_t\cap\ker\cC_t,
\qquad
\cM_k=\mathscr R_k\cap\ker\cC_k
\cap\cS_k^{-1}\cM_{k+1},
\quad s\le k<t,
\label{eq:backward-discard-recursion}
\end{equation}
where
$\cS_k^{-1}\cM_{k+1}
=\{h\in\cH_k:\cS_kh\in\cM_{k+1}\}$.
Then
\begin{equation}
\cM_k=\mathscr R_k\cap\ker\mathscr O_{t,k},
\quad s\le k\le t,
\qquad
\cS_k\cM_k\subseteq\cM_{k+1},
\quad s\le k<t.
\label{eq:finite-horizon-invisible-space}
\end{equation}

Set
$r_k=\operatorname{codim}_{\mathscr R_k}\cM_k$.
The quotient spaces
$\mathscr R_k/\cM_k$, $s\le k\le t$,
admit an exact linear state-space
realization that preserves all prescribed
task values $\cC_i h_i$ for every
initial history $h_k\in\mathscr R_k$,
every subsequent input sequence, and every $k\le i\le t$.
We call $\mathscr R_k/\cM_k$ the
reachable task quotient at time $k$
over $[k,t]$.

Moreover, if a linear map
$L_k:\mathscr R_k\to\F^{n_k}$
determines the zero-input future tasks,
in the sense that
$\mathscr O_{t,k}|_{\mathscr R_k}
=\Lambda_{t,k}L_k$
for some linear map
$\Lambda_{t,k}:\F^{n_k}
\to\F^{p_k+\cdots+p_t}$, then
\begin{equation}
n_k\ge
\operatorname{codim}_{\mathscr R_k}\cM_k
=r_k.
\label{eq:finite-horizon-dimension-lower-bound}
\end{equation}
The quotient realization attains
this lower bound. Hence $r_k$ is
the minimum dimension among linear
history states on $\mathscr R_k$
that determine all prescribed
zero-input future tasks.

Finally, let
$\mathscr G_{k,s}:
\F^{q_{k,s}}\to\mathscr R_k$
be any surjective linear map and define
$\mathscr H_{t,k,s}
=\mathscr O_{t,k}\mathscr G_{k,s}$.
Then
\begin{equation}
r_k
=
\rank\bigl(
\mathscr O_{t,k}|_{\mathscr R_k}
\bigr)
=
\rank\mathscr H_{t,k,s}.
\label{eq:reachable-hankel-rank}
\end{equation}
We call $r_k$ the future-task rank
at time $k$ over $[k,t]$.
It depends only on $\mathscr R_k$
and $\mathscr O_{t,k}$, and is
independent of the particular
surjective parameterization
$\mathscr G_{k,s}$.
\end{theorem}

\begin{proof}
We first establish
\eqref{eq:finite-horizon-invisible-space}
by backward induction.
At $k=t$, the identity follows from
$\mathscr O_{t,t}=\cC_t$.
Suppose that
$\cM_{k+1}=\mathscr R_{k+1}
\cap\ker\mathscr O_{t,k+1}$.
For every $h\in\mathscr R_k$,
the definition of the future-task
operator gives
$\mathscr O_{t,k}h=0$
if and only if $\cC_kh=0$ and
$\mathscr O_{t,k+1}\cS_kh=0$.
Since
$\cS_k\mathscr R_k
\subseteq\mathscr R_{k+1}$,
the induction hypothesis implies
\[
\mathscr R_k\cap\ker\mathscr O_{t,k}
=
\mathscr R_k\cap\ker\cC_k
\cap\cS_k^{-1}\cM_{k+1}
=\cM_k.
\]
The backward recursion also gives
$\cS_k\cM_k\subseteq\cM_{k+1}$,
which completes the induction.

Let
$q_k^{\mathrm{fh}}:
\mathscr R_k\to\mathscr R_k/\cM_k$
be the quotient map.
Since
$\cS_k\cM_k\subseteq\cM_{k+1}$
and $\cM_k\subseteq\ker\cC_k$,
the induced operators
\[
\widetilde{\cS}_k(q_k^{\mathrm{fh}}h)
=q_{k+1}^{\mathrm{fh}}\cS_kh,
\qquad
\widetilde{\cC}_k(q_k^{\mathrm{fh}}h)
=\cC_kh
\]
are well defined.
Indeed, two representatives of the
same quotient class differ by an
element of $\cM_k$, whose transported
image belongs to $\cM_{k+1}$
and whose current task value is zero.
Moreover,
$\operatorname{range}\cE_{k+1}
\subseteq\mathscr R_{k+1}$, so
$\widetilde{\cE}_{k+1}
=q_{k+1}^{\mathrm{fh}}\cE_{k+1}$
is well defined.

For every $h_k\in\mathscr R_k$
and entering input $u_{k+1}$,
\eqref{eq:history-migration} yields
\[
q_{k+1}^{\mathrm{fh}}h_{k+1}
=
\widetilde{\cS}_k
q_k^{\mathrm{fh}}h_k
+
\widetilde{\cE}_{k+1}u_{k+1},
\]
with
$\cC_kh_k
=\widetilde{\cC}_kq_k^{\mathrm{fh}}h_k$.
Because
$\cS_k\mathscr R_k+
\operatorname{range}\cE_{k+1}
=\mathscr R_{k+1}$,
the subsequent histories remain
in their corresponding reachable
subspaces for arbitrary entering inputs.

To verify task preservation explicitly,
consider two initial histories
$h,\hat h\in\mathscr R_k$
with $h-\hat h\in\cM_k$,
subject to the same subsequent input sequence.
Their difference at time $i\ge k$
is $\cS(i,k)(h-\hat h)$.
By
\eqref{eq:finite-horizon-invisible-space}
and the invariance of the subspaces
$\cM_j$,
this difference belongs to $\cM_i$
for every $k\le i\le t$.
Since $\cM_i\subseteq\ker\cC_i$,
the corresponding task values coincide.
Thus the quotient preserves all
prescribed tasks under arbitrary
subsequent inputs, even though
$\cM_k$ is characterized by the
zero-input future-task operator.

Choose linear isomorphisms
$U_k:\mathscr R_k/\cM_k\to\F^{r_k}$
and define
$x_k=U_kq_k^{\mathrm{fh}}h_k$.
For $s\le k<t$, set
\[
A_k=U_{k+1}\widetilde{\cS}_kU_k^{-1},
\qquad
B_{k+1}=U_{k+1}\widetilde{\cE}_{k+1}.
\]
Then
$x_{k+1}=A_kx_k+B_{k+1}u_{k+1}$
and
$\cC_kh_k
=\widetilde{\cC}_kU_k^{-1}x_k$.
Hence the quotient family provides
an exact linear state-space realization
of the prescribed tasks.
In particular, setting all subsequent
inputs to zero shows that $x_k$
determines $\mathscr O_{t,k}h_k$.

To establish minimality, suppose
$\mathscr O_{t,k}|_{\mathscr R_k}
=\Lambda_{t,k}L_k$.
For every $h\in\ker L_k$,
we have $\mathscr O_{t,k}h=0$.
Thus $\ker L_k\subseteq\cM_k$.
By rank--nullity,
\[
n_k\ge\rank L_k
=\operatorname{codim}_{\mathscr R_k}
\ker L_k
\ge
\operatorname{codim}_{\mathscr R_k}
\cM_k=r_k.
\]
This proves
\eqref{eq:finite-horizon-dimension-lower-bound}.
The constructed quotient state has
dimension $r_k$ and determines
all zero-input future tasks, so it
attains the lower bound.
Any linear realization preserving
the tasks under arbitrary subsequent
inputs must also preserve them when
those inputs vanish, and therefore
satisfies the same lower bound.

Finally,
$\ker(\mathscr O_{t,k}|_{\mathscr R_k})
=\cM_k$, and hence
\(\rank(
\mathscr O_{t,k}|_{\mathscr R_k})
=
\operatorname{codim}_{\mathscr R_k}
\cM_k=r_k.\)
Since $\mathscr G_{k,s}$ is
surjective onto $\mathscr R_k$,
\[
\operatorname{range}\mathscr H_{t,k,s}
=
\mathscr O_{t,k}
(\operatorname{range}\mathscr G_{k,s})
=
\operatorname{range}
(\mathscr O_{t,k}|_{\mathscr R_k}).
\]
Therefore
$\rank\mathscr H_{t,k,s}
=\rank(\mathscr O_{t,k}|_{\mathscr R_k})$,
which proves
\eqref{eq:reachable-hankel-rank}.
The rank depends only on
$\mathscr R_k$ and $\mathscr O_{t,k}$,
and is independent of the chosen
surjective parameterization.
\end{proof}

The reachable task quotient also admits
an explicit realization in orthonormal
coordinates. For each $s\le k\le t$,
let $Q_k:\F^{r_k}\to\cH_k$ be an
isometry with
$\operatorname{range}Q_k
=\mathscr R_k\cap\cM_k^\perp$.
Every $h_k\in\mathscr R_k$ admits
a unique decomposition
$h_k=Q_kx_k+m_k$, where
$m_k\in\cM_k$ and $x_k=Q_k^*h_k$.
Since $\cS_k\cM_k\subseteq\cM_{k+1}$,
the quotient coordinates satisfy
$x_{k+1}=A_kx_k+B_{k+1}u_{k+1}$, where
\begin{equation}
A_k=Q_{k+1}^*\cS_kQ_k,\qquad
B_{k+1}=Q_{k+1}^*\cE_{k+1}.
\label{eq:orthonormal-quotient-realization}
\end{equation}
The prescribed task values are recovered by
$\cC_kh_k=\cC_kQ_kx_k$.
In orthonormal coordinate bases,
$Q_k^*$ is represented by $Q_k^{\Hh}$.

In the finite-dimensional time-invariant setting,
the infinite-horizon reachable and task-invisible
subspaces admit finite-step characterizations.

\begin{corollary}
\label{cor:stationary-infinite-quotient}
Assume that $\cH_k=\cH_0$ for all $k$, with
$\dim\cH_0=d<\infty$, and that
$\cS\in\mathcal L(\cH_0,\cH_0)$,
$\cE\in\mathcal L(\F^m,\cH_0)$, and
$\cC\in\mathcal L(\cH_0,\F^p)$
are time invariant. Define
\begin{equation}
\mathscr R
=
\sum_{r=0}^{d-1}
\cS^r\operatorname{range}\cE,
\qquad
\cM
=
\mathscr R\cap
\bigcap_{r=0}^{d-1}\ker(\cC\cS^r).
\label{eq:stationary-reachable-quotient}
\end{equation}
Then
$\mathscr R
=\sum_{r=0}^{\infty}
\cS^r\operatorname{range}\cE$,
$\cM=\mathscr R\cap
\bigcap_{r=0}^{\infty}\ker(\cC\cS^r)$,
and $\cS\cM\subseteq\cM$.
Consequently, $\mathscr R/\cM$ admits
an exact time-invariant task-preserving
realization that determines
$\cC\cS^rh$ for every
$h\in\mathscr R$ and $r\ge0$.

Moreover, if a linear map
$L:\mathscr R\to\F^n$ satisfies
$\cC\cS^r|_{\mathscr R}=\Lambda_rL$
for some linear map
$\Lambda_r:\F^n\to\F^p$
at every $r\ge0$, then
$n\ge\dim(\mathscr R/\cM)$.
Thus $\mathscr R/\cM$ has the minimum
dimension among linear states preserving
all zero-input future task values on
$\mathscr R$. This dimension satisfies
\begin{equation}
\dim(\mathscr R/\cM)
=
\rank\mathscr H_d,
\qquad
\mathscr H_d
=
\bigl[\cC\cS^{i+j}\cE\bigr]_{i,j=0}^{d-1}.
\label{eq:stationary-hankel-rank}
\end{equation}
\end{corollary}

\begin{proof}
By the Cayley--Hamilton theorem,
every power $\cS^r$ with $r\ge d$
lies in the linear span of
$I_{\cH_0},\cS,\ldots,\cS^{d-1}$.
Consequently,
\[
\begin{aligned}
\mathscr R
&=\sum_{r=0}^{\infty}\cS^r\operatorname{range}\cE,
\qquad
\bigcap_{r=0}^{d-1}\ker(\cC\cS^r)
=\bigcap_{r=0}^{\infty}\ker(\cC\cS^r).
\end{aligned}
\]
These identities give the stated
infinite-horizon characterization
of $\cM$.

The representation of $\mathcal R$
shows that
$\mathcal S\mathcal R\subseteq\mathcal R$. For $h\in\cM$,
one has $\cS h\in\mathscr R$ and
$\cC\cS^r(\cS h)=0$ for every $r\ge0$.
Hence $\cS\cM\subseteq\cM$.

Let $q:\mathscr R\to\mathscr R/\cM$
be the quotient map. Since
$\cS\cM\subseteq\cM$ and
$\cM\subseteq\ker\cC$, the induced maps
$\widetilde{\cS}(qh)=q\cS h$
and
$\widetilde{\cC}(qh)=\cC h$
are well defined.
Moreover,
$\operatorname{range}\cE\subseteq\mathscr R$,
so $\widetilde{\cE}=q\cE$
is well defined.
For $h\in\mathscr R$ and
$u\in\F^m$, these maps satisfy
$q(\cS h+\cE u)
=\widetilde{\cS}(qh)+
\widetilde{\cE}u$
and
$\widetilde{\cC}(qh)=\cC h$.
Thus they define an exact
time-invariant task-preserving realization.
In particular,
$\widetilde{\cC}\widetilde{\cS}^{\,r}(qh)
=\cC\cS^rh$ for every $r\ge0$.

Suppose that
$\cC\cS^r|_{\mathscr R}=\Lambda_rL$
for every $r\ge0$.
If $h\in\ker L$, then
$\cC\cS^rh=0$ for every $r\ge0$,
so $h\in\cM$.
Therefore $\ker L\subseteq\cM$,
and rank--nullity gives
\[
n\ge\rank L
=
\operatorname{codim}_{\mathscr R}\ker L
\ge
\operatorname{codim}_{\mathscr R}\cM
=
\dim(\mathscr R/\cM).
\]
The quotient realization attains this
bound and therefore has minimum dimension
within the stated class.

Finally, define
$\mathscr O:\mathscr R\to\F^{pd}$
by
$\mathscr Oh=
\operatorname{col}(\cC h,\cC\cS h,
\ldots,\cC\cS^{d-1}h)$.
By \eqref{eq:stationary-reachable-quotient},
$\ker\mathscr O=\cM$, and hence
$\rank\mathscr O=\dim(\mathscr R/\cM)$.

Define
$\mathscr G:\F^{md}\to\mathscr R$
by
$\mathscr G=
[\cE,\cS\cE,\ldots,\cS^{d-1}\cE]$.
The definition of $\mathscr R$
implies that $\mathscr G$ is surjective,
so
$\operatorname{range}(\mathscr O\mathscr G)
=\operatorname{range}\mathscr O$.
Moreover, the $(i,j)$ block of
$\mathscr O\mathscr G$ is
$\cC\cS^{i+j}\cE$ for
$0\le i,j<d$.
Thus
$\mathscr O\mathscr G=\mathscr H_d$,
and consequently
$\rank\mathscr H_d
=\rank\mathscr O
=\dim(\mathscr R/\cM)$,
which proves
\eqref{eq:stationary-hankel-rank}.
\end{proof}

Unlike the general reachable future-task
matrix $\mathscr H_{t,k,s}$ in
Theorem~\ref{thm:finite-horizon-quotient-synthesis},
$\mathscr H_d$ is a block Hankel matrix
formed from the Markov parameters
$\cC\cS^r\cE$.
Its rank equals the classical minimum realization
order of the corresponding
input--output map under zero initial
conditions
\cite{Kalman1963,ho1966realization}.
The finite-step characterization follows
from the Cayley--Hamilton theorem and does
not extend without further assumptions
to general infinite-dimensional or
time-varying history dynamics.

\section{Minimum exact augmentation under time-varying history transport}\label{sec:timevaring}
Fix integers $s<t$ and assume that
$\cH_k$ is finite dimensional for
$s\le k\le t$. All minimality results
in this section concern exact linear
realizations over this fixed finite
horizon and on the full history spaces.
Throughout this section, a fixed basis
is chosen for each history space,
and all histories and operators are
represented in these bases.
With a slight abuse of notation,
we also use $h_k$ for the coordinate
vector of the history. Let $K_k$, $S_k$,
and $E_{k+1}$ denote the matrix
representations of $\cK_k$, $\cS_k$,
and $\cE_{k+1}$, respectively.
In these coordinates,
\eqref{eq:history-migration} becomes
\[
h_{k+1}=S_kh_k+E_{k+1}u_{k+1}.
\]

We seek linear history states
$\hat x_k=\hat K_kh_k$ that recover
the prescribed analysis $K_kh_k$
for every $h_k\in\cH_k$ and admit
an exact recurrence driven only by
the current state and the entering
input, for every such history and
every entering input.
Exact recovery requires
$K_k=F_k\hat K_k$ for some matrix $F_k$.
Moreover, the history dynamics give
\[
\hat x_{k+1}
=
\hat K_{k+1}S_kh_k
+\hat K_{k+1}E_{k+1}u_{k+1}.
\]
Thus an exact recurrence
$\hat x_{k+1}=A_k\hat x_k+B_{k+1}u_{k+1}$
for every history and input requires
$\hat K_{k+1}S_k=A_k\hat K_k$, with
$B_{k+1}=\hat K_{k+1}E_{k+1}$.
Conversely, these identities ensure
the required exact recurrence.

Accordingly, a matrix family
$\{\hat K_k\}_{k=s}^t$ is called
admissible if there exist matrices
$F_k$ and $A_k$ such that
\begin{equation}
K_k=F_k\hat K_k,\quad s\le k\le t,
\qquad
\hat K_{k+1}S_k=A_k\hat K_k,\quad s\le k<t.
\label{eq:augmentation-requirements}
\end{equation}
The resulting dimension bounds do not
extend to nonlinear state representations
or recurrences with independent access
to additional history information.

Equivalently,
\eqref{eq:augmentation-requirements}
can be expressed as
\[
\row K_k\subseteq\row\hat K_k,
\quad s\le k\le t,
\qquad
\row(\hat K_{k+1}S_k)
\subseteq\row\hat K_k,
\quad s\le k<t.
\]
Since these conditions depend only on
the row spaces of $\hat K_k$, we may
choose each $\hat K_k$ with linearly
independent rows when minimizing the
state dimension. Thus it suffices to
determine the smallest admissible
family of row spaces.

Define the backward invariant envelope
of the prescribed analyses as the
family of row spaces
$\{\cV_k\}_{k=s}^t$ given by
\begin{equation}
 \cV_t=\row K_t,
 \qquad
 \cV_k=\row K_k+\cV_{k+1}S_k,
 \quad s\le k<t,
 \label{eq:time-varying-invariant-envelope}
\end{equation}
where
$\cV_{k+1}S_k
=\{vS_k:v\in\cV_{k+1}\}$.
It is the smallest family of
row spaces containing $\row K_k$
and satisfying
$\cV_{k+1}S_k\subseteq\cV_k$
for $s\le k<t$.

\begin{theorem}
\label{thm:minimal-window-augmentation}
Let $\cV_k$ be defined by
\eqref{eq:time-varying-invariant-envelope}.
Then every admissible family
$\{\hat K_k\}_{k=s}^t$ satisfies
$\cV_k\subseteq\row\hat K_k$
for all $s\le k\le t$.
Conversely, choosing the rows of each
$\hat K_k$ to form a basis of $\cV_k$
produces an admissible family for all
$s\le k\le t$.
Consequently, over the fixed horizon
$[s,t]$ and within the admissible class
defined by~\eqref{eq:augmentation-requirements},
the minimum dimension of an exact
realization state $\hat x_k$ on the
full history spaces is $\dim\cV_k$.
The minimum number of additional
independent state directions beyond
$\row K_k$ is
$\dim\cV_k-\rank K_k$.
\end{theorem}

\begin{proof}
Let $\hat K_k$ be any admissible family.
At the terminal time,
$\cV_t=\row K_t\subseteq\row\hat K_t$.
Assume that
$\cV_{k+1}\subseteq\row\hat K_{k+1}$.
Then
\(\cV_{k+1}S_k
\subseteq
\row(\hat K_{k+1}S_k)
\subseteq
\row\hat K_k.\)
Since also
$\row K_k\subseteq\row\hat K_k$,
it follows that
\(\cV_k
=
\row K_k+\cV_{k+1}S_k
\subseteq
\row\hat K_k.\)
Backward induction therefore gives
$\cV_k\subseteq\row\hat K_k$
for all $s\le k\le t$.
Thus every admissible realization has
state dimension at least $\dim\cV_k$.

Conversely, choose each $\hat K_k$
with linearly independent rows spanning
$\cV_k$.
Since
$\row K_k\subseteq\cV_k=\row\hat K_k$,
there exists $F_k$ such that
$K_k=F_k\hat K_k$.
Moreover,
\[
\row(\hat K_{k+1}S_k)
=
\cV_{k+1}S_k
\subseteq
\cV_k
=
\row\hat K_k,
\]
so there exists $A_k$ satisfying
$\hat K_{k+1}S_k=A_k\hat K_k$.
Hence the family is admissible and has
exactly $\dim\cV_k$ state coordinates
at each time $k$.

Finally, since
$\row K_k\subseteq\cV_k$, a basis of
$\row K_k$ can be extended to a basis
of $\cV_k$ using exactly
$\dim\cV_k-\rank K_k$ additional vectors.
Because every admissible row space
contains $\cV_k$, fewer additional
independent directions cannot suffice.
This proves the claimed minimum
augmentation dimension within the
admissible class.
\end{proof}

Theorem~\ref{thm:finite-horizon-quotient-synthesis}
and Theorem~\ref{thm:minimal-window-augmentation}
address two distinct minimality problems.
The former concerns preservation of
future tasks on reachable histories,
whereas the latter concerns exact
recovery of prescribed analyses on
full history spaces within the
admissible class.

To relate these problems, set
$\cC_i=\cK_i$ for $k\le i\le t$,
with coordinate matrices $K_i$.
The definition of $\mathscr O_{t,k}$
and the backward recursion
\eqref{eq:time-varying-invariant-envelope}
then give
$\row\mathscr O_{t,k}=\cV_k$.
Consequently, the minimum state dimension
on the full history spaces within the
admissible class is
$d_k=\dim\cV_k=\rank\mathscr O_{t,k}$,
and the minimum number of additional
independent state directions is
$a_k=d_k-\rank K_k$.
The dimension $d_k$ allows the state
coordinates to be redesigned.
If all $N_k$ original analysis
coordinates must be retained as part
of the augmented state, the minimum
total number of state coordinates
is $N_k+a_k$.

In contrast, the minimum reachable
task dimension is
$r_k=\rank(\mathscr O_{t,k}|_{\mathscr R_k})
\le d_k$,
with equality when $\mathscr R_k=\cH_k$.
These dimensions concern a fixed
finite horizon and may depend on
the terminal time $t$.

\begin{corollary}
\label{cor:propagated-closure-defect-rank}
In the setting of
Theorem~\ref{thm:minimal-window-augmentation},
let $\Pi_j^{\mathrm{row}}$ be the orthogonal projector onto
$\operatorname{range}(K_j^{\Hh})$ with respect to
the standard inner product in the chosen history
coordinates, and define
$\Delta_j^{\mathrm{cl}}=K_{j+1}S_j(I-\Pi_j^{\mathrm{row}})$
for $s\le j<t$.
For $s\le k<t$, define the
propagated closure defect operator by
\(\mathsf D_{t,k}^{\mathrm{cl}}
=
\operatorname{col}\bigl(
\Delta_j^{\mathrm{cl}}S(j,k):
j=k,\ldots,t-1
\bigr),\)
where $S(j,k)=S_{j-1}\cdots S_k$ for $j>k$
and $S(k,k)=I$.
Then
\begin{align}
\cV_k
&=\row K_k+\row\mathsf D_{t,k}^{\mathrm{cl}},
\label{eq:propagated-defect-envelope}\\
\dim\cV_k-\operatorname{rank}K_k
&=\operatorname{rank}\bigl(
\mathsf D_{t,k}^{\mathrm{cl}}\big|_{\ker K_k}
\bigr).
\label{eq:propagated-defect-memory-rank}
\end{align}
For computation, let $\mathsf Q_k^{\mathrm{ker}}$
have orthonormal columns spanning $\ker K_k$
in the chosen history coordinates. Then
\(\operatorname{rank}\bigl(
\mathsf D_{t,k}^{\mathrm{cl}}\big|_{\ker K_k}
\bigr)
=
\operatorname{rank}\bigl(
\mathsf D_{t,k}^{\mathrm{cl}}\mathsf Q_k^{\mathrm{ker}}
\bigr),\)
independently of the choice of
the orthonormal basis.
Consequently,
\begin{equation}
0\le\dim\cV_k-\operatorname{rank}K_k
\le\min\left\{
\dim\ker K_k,\,
\sum_{j=k}^{t-1}\operatorname{rank}\Delta_j^{\mathrm{cl}}
\right\}.
\label{eq:propagated-defect-memory-bound}
\end{equation}
\end{corollary}
\begin{proof}
Expanding the backward recursion
\eqref{eq:time-varying-invariant-envelope}
gives
\(\cV_k
=
\sum_{i=k}^{t}
\row\bigl(K_iS(i,k)\bigr).\)
We show that this row space is generated
by the current analysis and the propagated
one-step closure defects.

Since
$\row\Pi_j^{\mathrm{row}}=\row K_j$,
there exists a matrix $\Gamma_j$ such that
$K_{j+1}S_j\Pi_j^{\mathrm{row}}
=\Gamma_jK_j$.
Consequently,
$K_{j+1}S_j
=\Gamma_jK_j+\Delta_j^{\mathrm{cl}}$.
Right multiplication by $S(j,k)$ gives
\[
K_{j+1}S(j+1,k)
=
\Gamma_jK_jS(j,k)
+
\Delta_j^{\mathrm{cl}}S(j,k).
\]
Hence, by elementary row operations on block matrices,
\[
\row\operatorname{col}\bigl(
K_jS(j,k),K_{j+1}S(j+1,k)
\bigr)
=
\row\operatorname{col}\bigl(
K_jS(j,k),\Delta_j^{\mathrm{cl}}S(j,k)
\bigr).
\]
Applying this identity successively
for $j=k,\ldots,t-1$ yields
\[
\sum_{i=k}^{t}
\row\bigl(K_iS(i,k)\bigr)
=
\row K_k+
\sum_{j=k}^{t-1}
\row\bigl(
\Delta_j^{\mathrm{cl}}S(j,k)
\bigr).
\]
By the definition of
$\mathsf D_{t,k}^{\mathrm{cl}}$,
this proves
\eqref{eq:propagated-defect-envelope}.
It follows that
$\dim\cV_k
=\rank\operatorname{col}
(K_k,\mathsf D_{t,k}^{\mathrm{cl}})$.
Applying Lemma~\ref{lem:rank-increment}
with $K=K_k$ and
$L=\mathsf D_{t,k}^{\mathrm{cl}}$
gives
\eqref{eq:propagated-defect-memory-rank}.

Since the columns of
$\mathsf Q_k^{\mathrm{ker}}$
span $\ker K_k$, we have
\(\operatorname{range}\bigl(
\mathsf D_{t,k}^{\mathrm{cl}}
\mathsf Q_k^{\mathrm{ker}}
\bigr)
=
\mathsf D_{t,k}^{\mathrm{cl}}(\ker K_k).\)
Hence the stated computational rank
identity follows. Any two orthonormal
bases of $\ker K_k$ are related by a
unitary transformation, so the rank
is independent of the basis choice.
Finally,
\(\rank\bigl(
\mathsf D_{t,k}^{\mathrm{cl}}
\big|_{\ker K_k}
\bigr)
\le
\min\left\{
\dim\ker K_k,\,
\rank\mathsf D_{t,k}^{\mathrm{cl}}
\right\}.\)
By rank subadditivity and
$\rank(\Delta_j^{\mathrm{cl}}S(j,k))
\le\rank\Delta_j^{\mathrm{cl}}$,
we obtain
\(\rank\mathsf D_{t,k}^{\mathrm{cl}}
\le
\sum_{j=k}^{t-1}
\rank\Delta_j^{\mathrm{cl}}.\)
This proves
\eqref{eq:propagated-defect-memory-bound}.
\end{proof}

The restricted rank in
\eqref{eq:propagated-defect-memory-rank}
counts the history directions invisible
to the current analysis that become
relevant under future transport.
Propagated defects may contain
overlapping information, so summing
their individual ranks can overestimate
the minimum augmentation dimension.
The following example shows that the
bound in
\eqref{eq:propagated-defect-memory-bound}
can be strict.
\begin{example}
    Let $s=0$, $t=4$, $\cH_j=\R^6$, and
\[
K_j=K=[1,0_{1\times5}],
\qquad
S_j=
\begin{bmatrix}
1&v_j&0_{1\times3}\\
0_{2\times1}&R&0_{2\times3}\\
0_{3\times1}&0_{3\times2}&I_3
\end{bmatrix},
\]
where
$R=\begin{bmatrix}1&1\\0&1\end{bmatrix}$
and $v_j\in\R^{1\times2}$ is defined by
$v_0=v_2=[1,0]$ and
$v_1=v_3=[0,1]$.
Since
$v_1+v_0R\ne v_0+v_1R$,
the matrices $S_0$ and $S_1$
do not commute.
As $\Pi_j^{\mathrm{row}}=K^{\T}K$,
the one-step defects and their
propagations satisfy
\[
\Delta_j^{\mathrm{cl}}
=[0,v_j,0_{1\times3}],
\qquad
\Delta_j^{\mathrm{cl}}S(j,0)
=[0,v_jR^j,0_{1\times3}].
\]
For $j=0,1,2,3$, the vectors
$v_jR^j$ are
$[1,0]$, $[0,1]$, $[1,2]$, and
$[0,1]$, respectively.
They span a two-dimensional space.
Since $\ker K=\{h\in\R^6:h_1=0\}$
and the propagated defect rows
involve only the second and third
coordinates, their restriction to
$\ker K$ retains rank two.
Consequently,
\[
\rank\bigl(
\mathsf D_{4,0}^{\mathrm{cl}}
\big|_{\ker K}
\bigr)
=2<4
=\sum_{j=0}^{3}\rank\Delta_j^{\mathrm{cl}}.
\]
By
\eqref{eq:propagated-defect-memory-rank},
the minimum additional dimension
at time zero is two, giving
$\dim\cV_0=3$.
This minimum is attained by
$\hat K_j=[I_3,0_{3\times3}]$, since
\[
\hat K_{j+1}S_j=\hat A_j\hat K_j,
\qquad
\hat A_j=
\begin{bmatrix}
1&v_j\\
0_{2\times1}&R
\end{bmatrix},
\qquad
K=[1,0,0]\hat K_j.
\]
Thus four nonzero one-step defects
require only two additional state
directions for exact recovery.
\end{example}

The projector $\Pi_j^{\mathrm{row}}$
is defined using the standard inner product
in the chosen history coordinates.
Consequently, $\Pi_j^{\mathrm{row}}$,
$\Delta_j^{\mathrm{cl}}$, and their matrix
norms generally depend on the coordinate
choice. In contrast, the restricted rank
in \eqref{eq:propagated-defect-memory-rank}
and the minimum augmentation dimension
are invariant under invertible changes
of history coordinates.
If the chosen coordinates are orthonormal
with respect to the inner products on
the history spaces and $K_jK_j^{\Hh}=I$,
then $\Delta_j^{\mathrm{cl}}$ represents
the cross-projection defect
$\Delta_j^{\mathrm{cp}}$ in
\eqref{eq:cross-projection-defect}.

For a fixed $d$-dimensional history space
with time-independent transport $S$
and prescribed analysis $K$, the backward
recursion gives
$\cV_k=\sum_{r=0}^{t-k}\row(KS^r)$.
By the Cayley--Hamilton theorem, this
space stabilizes when $t-k\ge d-1$,
with
\begin{equation}
\cV(K;S)
=
\sum_{r=0}^{d-1}\row(KS^r).
\label{eq:general-invariant-envelope}
\end{equation}
The space $\cV(K;S)$ is the row space
of the finite-step observability matrix
associated with $(S,K)$
\cite{Kalman1963,zhou}.
It is the smallest row space containing
$\row K$ and invariant under right
multiplication by $S$.

\section{Structured exact realizations for countable histories}
\label{sec:structured-realizations}

The preceding sections characterize
exact closure and finite-horizon
minimality for prescribed history
analyses. We now construct a class
of analyses admitting exact closure
on geometrically weighted countable
histories using possibly noncommuting
matrix modes, and establish state
transition bounds and sufficient
conditions for uniform Gram bounds.
For a stationary Meixner specialization,
restriction to finite causal histories
allows the minimum augmentation theory
to be applied to both the full analysis
of retained coefficients and individual
prescribed coefficients.

\subsection{Matrix mode transport}
\label{subsec:matrix-mode-dho}

We use a geometrically weighted Hilbert history
space and prescribe a family of analysis kernels
through matrix mode transport.

Fix an integer $N\ge1$ and
$0<\sigma<1$. Here $N$ is the
fixed number of analysis coordinates,
corresponding to $N_k=N$ in the
general framework.
Define the probability weight
$\nu(\ell)=(1-\sigma)\sigma^\ell$
on $\N$.
Let
\[
 \ell^2(\nu)
 =
 \left\{
 f:\N\to\C:
 \sum_{\ell=0}^{\infty}\nu(\ell)|f(\ell)|^2<\infty
 \right\},
\]
equipped with the inner product
$\langle f,g\rangle_\nu
=
\sum_{\ell=0}^{\infty}
\nu(\ell)f(\ell)\overline{g(\ell)}$.

Let $b_k\in\C^N$ and $R_k\in\C^{N\times N}$
for every $k\in\mathbb Z$.
The bilateral indexing ensures that
$g_k(\ell)$ is well defined for every
$k\ge0$ and $\ell\ge0$, even though
the input sequence is causal.
For each $k\in\mathbb Z$, define
a vector-valued sequence
$g_k:\N\to\C^N$ recursively by
\begin{equation}
 g_k(0)=b_k,
 \qquad
 g_{k+1}(\ell+1)=R_{k+1}g_k(\ell),
 \qquad
 \ell\ge0.
 \label{eq:matrix-mode-transport}
\end{equation}
The matrices $R_k$ are not assumed to commute.
For $k\ge s\ge0$, define
$\mathscr P(k,s)=R_kR_{k-1}\cdots R_{s+1}$
when $k>s$, and set
$\mathscr P(s,s)=I_N$.
These ordered matrix products are related
to matrix product growth in switching
systems \cite{jungers2009jsr}.

Assume that, for every $k\ge0$,
the weighted Gram matrix
\begin{equation}
 G_k
 =
 \sum_{\ell=0}^{\infty}
 \nu(\ell)g_k(\ell)g_k^{\Hh}(\ell)
 \label{eq:matrix-mode-gram}
\end{equation}
is well defined and positive definite.
For $k\ge0$, define
$\Theta_k=G_k^{-1/2}$,
$\phi_k(\ell)=\Theta_k g_k(\ell)$, and write
\[
\phi_k(\ell)
=
\operatorname{col}\bigl(
\phi_{k,0}(\ell),\ldots,\phi_{k,N-1}(\ell)
\bigr).
\]
Then
$\sum_{\ell\ge0}\nu(\ell)
\phi_k(\ell)\phi_k^{\Hh}(\ell)=I_N$,
so the component sequences
$\{\phi_{k,n}\}_{n=0}^{N-1}$
form an orthonormal family in
$\ell^2(\nu)$.
Define
$\psi_{k,n}(\ell)=\overline{\phi_{k,n}(\ell)}$
for $0\le n<N$ and $\ell\ge0$,
where conjugation is pointwise.
Since complex conjugation preserves
orthonormality,
$\{\psi_{k,n}\}_{n=0}^{N-1}$
is also an orthonormal family in
$\ell^2(\nu)$.

For a scalar sequence
$\{u_k\}_{k\ge0}\subset\C$, set
$u_j=0$ for $j<0$ and define the lag
history by $h_k(\ell)=u_{k-\ell}$.
Then $h_k(\ell)=0$ for $\ell>k$, so
$h_k\in\ell^2(\nu)$.
For this history model, take
$\cH_k=\ell^2(\nu)$ and define
the shift operator
$\cS:\ell^2(\nu)\to\ell^2(\nu)$ by
$(\cS h)(0)=0$ and
$(\cS h)(\ell)=h(\ell-1)$ for $\ell\ge1$,
and define the insertion operator
$\cE:\C\to\ell^2(\nu)$ by
$(\cE u)(0)=u$ and
$(\cE u)(\ell)=0$ for $\ell\ge1$.
Then
$h_{k+1}=\cS h_k+\cE u_{k+1}$.
Moreover,
$\norm{\cS h}_{\ell^2(\nu)}^2
=\sigma\norm{h}_{\ell^2(\nu)}^2$,
so $\cS$ is bounded with
$\norm{\cS}=\sqrt{\sigma}$.

Let $\cK_k:\ell^2(\nu)\to\C^N$
be defined by
$[\cK_kh]_n=\langle h,\psi_{k,n}\rangle_\nu$,
$0\le n<N$.
Since the inner product is linear
in its first argument and
$\psi_{k,n}=\overline{\phi_{k,n}}$, we have
\[
[\cK_kh]_n
=\sum_{\ell=0}^{\infty}
\nu(\ell)h(\ell)\overline{\psi_{k,n}(\ell)}
=\sum_{\ell=0}^{\infty}
\nu(\ell)\phi_{k,n}(\ell)h(\ell).
\]
By orthonormality,
$\cK_k^*x=\sum_{n=0}^{N-1}x_n\psi_{k,n}$
and $\cK_k\cK_k^*=I_N$.
Thus $\cK_k$ is surjective.
The corresponding DHO state is
\begin{equation}
 x_k
 =
 \cK_kh_k
 =
 \sum_{\ell=0}^{\infty}
 \nu(\ell)\phi_k(\ell)u_{k-\ell}.
 \label{eq:matrix-mode-state}
\end{equation}

\begin{theorem}
\label{thm:matrix-mode-dho}
The state in \eqref{eq:matrix-mode-state}
satisfies, for $k\ge0$, the exact
DLTV recurrence
\begin{equation}
 x_{k+1}
 =
 A_kx_k+B_{k+1}u_{k+1},
 \qquad
 A_k
 =
 \sigma \Theta_{k+1}R_{k+1}\Theta_k^{-1},
 \qquad
 B_{k+1}
 =
 (1-\sigma)\Theta_{k+1}b_{k+1}.
 \label{eq:matrix-mode-matrices}
\end{equation}
For the state transition $\cF(k,s)$
defined in \eqref{eq:ssm-kernel},
we have, for $k\ge s\ge0$,
\begin{equation}
 \cF(k,s)
 =
 \sigma^{k-s}
 \Theta_k\mathscr P(k,s)\Theta_s^{-1}.
 \label{eq:matrix-mode-transition}
\end{equation}
Moreover, the state transitions satisfy
\[
\norm{\cF(k,s)}
\le \sigma^{(k-s)/2},
\qquad k\ge s\ge0.
\]
If, in addition, there exist constants
$0<\underline g\le\bar g$ such that
$\underline g I_N\preceq G_i
\preceq\bar g I_N$ for every $i\ge0$,
then
\begin{equation}
 \sqrt{\frac{\underline g}{\bar g}}\,
 \sigma^{k-s}
 \sigma_{\min}\!\bigl(\mathscr P(k,s)\bigr)
 \le
 \sigma_{\min}\!\bigl(\cF(k,s)\bigr)
 \le
 \norm{\cF(k,s)}
 \le
 \sqrt{\frac{\bar g}{\underline g}}\,
 \sigma^{k-s}
 \norm{\mathscr P(k,s)}.
 \label{eq:matrix-mode-transition-bounds}
\end{equation}
\end{theorem}

\begin{proof}
Separating the term $\ell=0$ in
\eqref{eq:matrix-mode-state}
at time $k+1$ gives
\[
 x_{k+1}
 =
 (1-\sigma)\Theta_{k+1}b_{k+1}u_{k+1}
 +
 \sum_{\ell=1}^{\infty}
 (1-\sigma)\sigma^\ell
 \Theta_{k+1}g_{k+1}(\ell)u_{k+1-\ell}.
\]
By \eqref{eq:matrix-mode-transport},
$g_{k+1}(\ell)=R_{k+1}g_k(\ell-1)$
for $\ell\ge1$.
Setting $q=\ell-1$, the second term becomes
$\sigma\Theta_{k+1}R_{k+1}
\Theta_k^{-1}x_k$.
This proves \eqref{eq:matrix-mode-matrices}.
The same transport identity gives
$\cK_{k+1}\cS=A_k\cK_k$ on
$\ell^2(\nu)$, while
$B_{k+1}=\cK_{k+1}\cE$.
Thus the operator closure condition
\eqref{eq:operator-closure} holds.

For $k>s$, \eqref{eq:ssm-kernel} and
\eqref{eq:matrix-mode-matrices} give
\[
\cF(k,s)
=
A_{k-1}\cdots A_s
=
\sigma^{k-s}
\Theta_kR_kR_{k-1}\cdots
R_{s+1}\Theta_s^{-1},
\]
which is \eqref{eq:matrix-mode-transition}.
For $k=s$, the same identity follows from
$\cF(s,s)=\mathscr P(s,s)=I_N$.

Since $\cK_{i+1}\cS=A_i\cK_i$
and $\cK_i\cK_i^*=I_N$, we have
$A_i=\cK_{i+1}\cS\cK_i^*$.
Thus
$\norm{A_i}\le\norm{\cS}=\sqrt{\sigma}$,
which yields
$\norm{\cF(k,s)}
\le\sigma^{(k-s)/2}$
for every $k\ge s\ge0$.

Since $\Theta_i=G_i^{-1/2}$ and
$\Theta_i^{-1}=G_i^{1/2}$, the
assumed uniform Gram bounds imply
$\norm{\Theta_i}\le1/\sqrt{\underline g}$,
$\sigma_{\min}(\Theta_i)\ge1/\sqrt{\bar g}$,
$\norm{\Theta_i^{-1}}\le\sqrt{\bar g}$,
and
$\sigma_{\min}(\Theta_i^{-1})
\ge\sqrt{\underline g}$.
Therefore \eqref{eq:matrix-mode-transition}
gives
\[
 \sigma_{\min}\!\bigl(\cF(k,s)\bigr)
 \ge
 \sqrt{\frac{\underline g}{\bar g}}\,
 \sigma^{k-s}
 \sigma_{\min}\!\bigl(\mathscr P(k,s)\bigr),
 \qquad
 \norm{\cF(k,s)}
 \le
 \sqrt{\frac{\bar g}{\underline g}}\,
 \sigma^{k-s}
 \norm{\mathscr P(k,s)}.
\]
Together with
$\sigma_{\min}(\cF(k,s))
\le\norm{\cF(k,s)}$,
this proves
\eqref{eq:matrix-mode-transition-bounds}.
\end{proof}

Exact closure and a uniform exponential
decay estimate for state transitions hold
without uniform Gram bounds.
Under additional uniform Gram bounds,
\eqref{eq:matrix-mode-transition-bounds}
relates the smallest singular value
and operator norm of the state transition
to those of the ordered mode product.
In \eqref{eq:matrix-mode-transition},
the factor $\sigma^{k-s}$ arises from
the geometric weight, while
$\mathscr P(k,s)$ retains the ordered,
possibly noncommuting mode transport.

\begin{theorem}
\label{thm:finite-gram-certificate}
Assume that there exist constants $r>0$
and $\beta\ge0$ such that $\norm{R_i}\le r$
and $\norm{b_i}\le\beta$ for every
$i\in\mathbb Z$, with $\sigma r^2<1$.
For an integer $L\ge N$, define
\begin{equation}
W_{k,L}
=
(1-\sigma)\sum_{\ell=0}^{L-1}
\sigma^\ell g_k(\ell)g_k^{\Hh}(\ell).
\label{eq:finite-mode-gramian}
\end{equation}
If there exists $\alpha>0$ such that
$W_{k,L}\succeq\alpha I_N$ for every $k\ge0$,
then
\begin{equation}
\alpha I_N
\preceq G_k
\preceq
\frac{(1-\sigma)\beta^2}
{1-\sigma r^2}I_N,
\qquad k\ge0.
\label{eq:finite-gram-certificate}
\end{equation}
Moreover, the state transition matrix
satisfies
\begin{equation}
\norm{\cF(k,s)}
\le
\sqrt{\frac{(1-\sigma)\beta^2}
{\alpha(1-\sigma r^2)}}
(\sigma r)^{k-s},
\qquad k\ge s\ge0.
\label{eq:finite-gram-transition-decay}
\end{equation}
In particular, the decay rate
$\sigma r$ is strictly smaller than
the general rate $\sqrt{\sigma}$.
\end{theorem}

\begin{proof}
From \eqref{eq:matrix-mode-transport},
$g_k(\ell)
=R_kR_{k-1}\cdots R_{k-\ell+1}b_{k-\ell}$
for $\ell\ge1$.
Hence $\norm{g_k(\ell)}\le r^\ell\beta$
for every $k\ge0$ and $\ell\ge0$.
Since every summand in
\eqref{eq:matrix-mode-gram}
is positive semidefinite,
$G_k\succeq W_{k,L}\succeq\alpha I_N$.
For any $z\in\C^N$ with $\norm{z}=1$,
\[
z^{\Hh}G_kz
\le
(1-\sigma)\beta^2
\sum_{\ell=0}^{\infty}
(\sigma r^2)^\ell
=
\frac{(1-\sigma)\beta^2}
{1-\sigma r^2}.
\]
Taking the supremum over unit vectors
gives \eqref{eq:finite-gram-certificate}.

Furthermore, the uniform bound
$\norm{R_i}\le r$ implies
$\norm{\mathscr P(k,s)}\le r^{k-s}$.
Combining this inequality with
\eqref{eq:matrix-mode-transition-bounds}
and \eqref{eq:finite-gram-certificate}
gives \eqref{eq:finite-gram-transition-decay}.
Since $\sigma r^2<1$ and $0<\sigma<1$,
we have $0<\sigma r<\sqrt{\sigma}<1$,
which completes the proof.
\end{proof}

The finite-lag lower bound requires
uniformity over all $k\ge0$ and cannot
be inferred from finitely many time
indices in general nonperiodic systems.
The conditions are sufficient but not
necessary for uniform Gram bounds.
Moreover, the Euclidean condition
$\sigma r^2<1$ may fail even when
a suitable common metric satisfies
$\sigma r_P^2<1$.

\begin{proposition}
\label{prop:metric-gram-certificate}
Let $P=P^{\Hh}\succ0$, and let $L\ge N$ be an integer.
Assume that there exist constants $r_P,\beta_P\ge0$ such that,
for every $i\in\mathbb Z$,
$R_i^{\Hh}PR_i\preceq r_P^2P$ and
$b_i^{\Hh}Pb_i\le\beta_P^2$, with $\sigma r_P^2<1$.
If there exists $\alpha>0$ such that
$W_{k,L}\succeq\alpha P^{-1}$ for every $k\ge0$, then
\begin{equation}
 \alpha P^{-1}
\preceq
 G_k
\preceq
 \frac{(1-\sigma)\beta_P^2}
 {1-\sigma r_P^2}
 P^{-1},
 \qquad k\ge0.
 \label{eq:metric-gram-bounds}
\end{equation}
\end{proposition}

\begin{proof}
For $\ell\ge1$,
$g_k(\ell)
=
R_kR_{k-1}\cdots R_{k-\ell+1}b_{k-\ell}$.
Repeated application of
$R_i^{\Hh}PR_i\preceq r_P^2P$ gives
$g_k(\ell)^{\Hh}Pg_k(\ell)
\le\beta_P^2r_P^{2\ell}$
for every $k\ge0$ and $\ell\ge0$.

For any $z\in\C^N$, the Cauchy--Schwarz inequality gives
\[
 \abs{z^{\Hh}g_k(\ell)}^2
 =
 \abs{
 (P^{-1/2}z)^{\Hh}
 (P^{1/2}g_k(\ell))
 }^2
 \le
 (z^{\Hh}P^{-1}z)
 (g_k(\ell)^{\Hh}Pg_k(\ell)).
\]
Hence
\[
 z^{\Hh}G_kz
 \le
 (1-\sigma)\beta_P^2
 z^{\Hh}P^{-1}z
 \sum_{\ell=0}^{\infty}
 (\sigma r_P^2)^\ell
 =
 \frac{(1-\sigma)\beta_P^2}
 {1-\sigma r_P^2}
 z^{\Hh}P^{-1}z.
\]
Since this holds for every $z\in\C^N$, it gives the upper inequality
in \eqref{eq:metric-gram-bounds}.
The lower inequality follows from
$G_k\succeq W_{k,L}\succeq\alpha P^{-1}$.
\end{proof}

In particular, \eqref{eq:metric-gram-bounds}
implies the uniform Gram hypothesis in
Theorem~\ref{thm:matrix-mode-dho} with
$\underline g=\alpha/\lambda_{\max}(P)$ and
$\bar g=(1-\sigma)\beta_P^2/
[(1-\sigma r_P^2)\lambda_{\min}(P)]$.

If $R_i$ and $b_i$ have a common
period $T$, then $W_{k+T,L}=W_{k,L}$.
Consequently, for fixed $P\succ0$
and $L\ge N$, there exists
$\alpha>0$ such that
$W_{k,L}\succeq\alpha P^{-1}$
for all $k\ge0$ if and only if
$W_{p,L}\succ0$ for $0\le p<T$.
Thus the uniform lower bound can
be established by checking finitely
many phases in the periodic case.

\subsection{Meixner polynomial realizations}
\label{subsec:meixner}
We next specialize the matrix mode framework to
a stationary family of discrete orthogonal
polynomials. Fix an order $N\ge1$, representing
the number of retained Meixner
coefficients, and retain the geometric
parameter $0<\sigma<1$.
The geometric lag weight
$\nu(\ell)=(1-\sigma)\sigma^\ell$,
$\ell\in\N$, remains unchanged.
This weight and its Meixner orthogonal basis
are independent of time.
As before, extend the input by $u_j=0$
for $j<0$.
For $n,\ell\in\N$, define the normalized
Meixner polynomials by
\begin{equation}
 m_n(\ell)
 =
 \sigma^{n/2}
 {}_2F_1\!\left(
 \begin{matrix}
 -n,-\ell\\
 1
 \end{matrix};
 1-\sigma^{-1}
 \right).
 \label{eq:meixner-basis}
\end{equation}
Here $(a)_0=1$ and $(a)_r=a(a+1)\cdots(a+r-1)$
for $r\ge1$. The terminating hypergeometric series is
\[
 {}_2F_1\!\left(
 \begin{matrix}
 -n,-\ell\\
 1
 \end{matrix};z
 \right)
 =
 \sum_{r=0}^{n}
 \frac{(-n)_r(-\ell)_r}{(1)_r}
 \frac{z^r}{r!}.
\]
These polynomials are orthonormal under the geometric weight
\cite{ismail2005classical,koekoek2010hypergeometric}, that is,
\[
 \sum_{\ell=0}^{\infty}
 (1-\sigma)\sigma^\ell
 m_n(\ell)m_q(\ell)
 =
 \delta_{nq}.
\]
Here $\delta_{nq}$ denotes the Kronecker delta.

For $0\le n<N$, define
\begin{equation}
 x_{k,n}
 =
 \sum_{\ell=0}^{\infty}
 (1-\sigma)\sigma^\ell
 m_n(\ell)u_{k-\ell},
 \label{eq:meixner-state}
\end{equation}
and let
$x_k=\operatorname{col}(x_{k,0},\ldots,x_{k,N-1})$.
Because $u_j=0$ for $j<0$, the sum in
\eqref{eq:meixner-state} contains only the terms
$\ell=0,\ldots,k$ at time $k$.

\begin{theorem}
\label{thm:meixner-realization}
The state in \eqref{eq:meixner-state} satisfies the exact
time-invariant DHO recurrence
\(x_{k+1}=A^{\mathrm M}x_k+B^{\mathrm M}u_{k+1}\),
where
\begin{equation}
 [A^{\mathrm M}]_{n,m}
 =
 \begin{cases}
 \sigma, & n=m,\\
 -(1-\sigma)\sigma^{(n-m)/2}, & n>m,\\
 0, & n<m,
 \end{cases}
 \qquad
 [B^{\mathrm M}]_n
 =
 (1-\sigma)\sigma^{n/2},
 \label{eq:meixner-matrices}
\end{equation}
for $0\le m,n<N$.
\end{theorem}
\begin{proof}
Let $M_n(\ell)$ denote the Meixner
polynomial in \cite{xu2024meixner}
with shape parameter one and
geometric parameter $\sigma$.
Under the normalization in
\eqref{eq:meixner-basis},
$m_n(\ell)=\sigma^{n/2}M_n(\ell)/n!$.
By \cite[Theorem~3.1]{xu2024meixner},
\[
M_n(\ell+1)-M_n(\ell)
=
(1-\sigma^{-1})
\sum_{j=0}^{n-1}
\frac{n!}{j!}M_j(\ell).
\]
Substitution of the normalization
gives
\begin{equation}
\sigma m_n(\ell+1)
=
\sigma m_n(\ell)
-(1-\sigma)
\sum_{j=0}^{n-1}
\sigma^{(n-j)/2}m_j(\ell).
\label{eq:meixner-connection}
\end{equation}

Let
$g_k(\ell)=\operatorname{col}
(m_0(\ell),\ldots,m_{N-1}(\ell))$,
independently of $k$.
By Meixner orthonormality,
the matrix mode construction gives
$G_k=I_N$ and $\Theta_k=I_N$.
Equation~\eqref{eq:meixner-connection}
implies
$g_{k+1}(\ell+1)=R^{\mathrm M}g_k(\ell)$,
where $R^{\mathrm M}=\sigma^{-1}A^{\mathrm M}$
and the entries of $A^{\mathrm M}$
are given by \eqref{eq:meixner-matrices}.
Moreover, $b_k=g_k(0)$ and
$[b_k]_n=m_n(0)=\sigma^{n/2}$.
Theorem~\ref{thm:matrix-mode-dho}
therefore gives
$A_k=\sigma R^{\mathrm M}=A^{\mathrm M}$
and
$B_{k+1}=(1-\sigma)b_{k+1}=B^{\mathrm M}$,
which proves the exact recurrence.
\end{proof}

Thus the Meixner realization is a
stationary specialization of the general
matrix mode construction, with $G_k=I_N$.
Consequently, uniform Gram bounds hold
directly in this case.

The difference identity in
\eqref{eq:meixner-connection}
is a normalized form of the known
Meixner relation
\cite[Theorem~3.1]{xu2024meixner}.
Recursive state-space realizations
of related Meixner-like filters
for transfer function models have
also been developed in
\cite{maraoui2017meixner}.
The distinction examined here is
the minimum exact augmentation of
the full Meixner analysis versus
a single prescribed coefficient
over finite causal histories, as
shown in Examples~\ref{example}
and~\ref{ex:meixner-selected}.

Since $A^{\mathrm M}$ is lower triangular
with diagonal entries $\sigma$, its
spectral radius is $\sigma<1$.
For $N>1$, however, the matrix is
nonnormal, and its eigenvalues alone
do not determine the norms of its powers.
We next derive an explicit matrix power
formula and an exponential bound for
fixed $N$.

\begin{proposition}
\label{prop:meixner-powers}
Let
$S_\sigma=\operatorname{diag}(1,\sigma^{1/2},\ldots,\sigma^{(N-1)/2})$
and let $J$ be the strictly lower triangular matrix whose entries
below the diagonal are all equal to one.
Then, for every $k\in\N$,
\begin{equation}
 (A^{\mathrm M})^k
 =
 S_\sigma
 \sum_{r=0}^{\min\{k,N-1\}}
 \binom{k}{r}
 \sigma^{k-r}
 [-(1-\sigma)J]^r
 S_\sigma^{-1}.
 \label{eq:meixner-matrix-power}
\end{equation}
Consequently, for every $\hat\sigma\in(\sigma,1)$, there exists a
constant $M_{N,\sigma,\hat\sigma}>0$ such that
\[
 \norm{(A^{\mathrm M})^k}
 \le
 M_{N,\sigma,\hat\sigma}\hat\sigma^k,
 \qquad
 k\in\N.
\]
\end{proposition}

\begin{proof}
By \eqref{eq:meixner-matrices},
$S_\sigma^{-1}A^{\mathrm M}S_\sigma
=\sigma I_N-(1-\sigma)J$.
Since $J$ is strictly lower triangular,
$J^N=0$.
The binomial theorem and similarity
transformation give
\eqref{eq:meixner-matrix-power}.

Fix $\hat\sigma\in(\sigma,1)$.
For each fixed $0\le r<N$ and $k\ge r$,
we have
$\binom{k}{r}\le (1+k)^r/r!$.
Since $\sigma/\hat\sigma<1$, the constant
\[
C_r:=
\frac{1}{r!}
\sup_{k\in\N}
(1+k)^r(\sigma/\hat\sigma)^k
\]
is finite. Hence
\[
\binom{k}{r}\sigma^{k-r}
\le C_r\sigma^{-r}\hat\sigma^k,
\qquad k\ge r.
\]
Applying this estimate to
\eqref{eq:meixner-matrix-power} yields
\[
\norm{(A^{\mathrm M})^k}
\le
\norm{S_\sigma}\norm{S_\sigma^{-1}}
\sum_{r=0}^{\min\{k,N-1\}}
C_r\sigma^{-r}(1-\sigma)^r
\norm{J^r}\hat\sigma^k\le M_{N,\sigma,\hat\sigma}\hat\sigma^k,
\]
where
\(M_{N,\sigma,\hat\sigma}
:=
\norm{S_\sigma}\norm{S_\sigma^{-1}}
\sum_{r=0}^{N-1}
C_r\sigma^{-r}(1-\sigma)^r\norm{J^r}.\)
Since $N$ is fixed and each $C_r$ is
finite, $M_{N,\sigma,\hat\sigma}$ is finite
and independent of $k$.
Thus the estimate holds for all
$k\in\N$.
\end{proof}

The structure of $A^{\mathrm M}$ permits
a state update in $O(N)$ arithmetic
operations using diagonal scaling
and prefix sums.

\begin{proposition}
\label{prop:meixner-structured-scaling}
With $S_\sigma$ and $J$ as in
Proposition~\ref{prop:meixner-powers}, define
$\widetilde x_k=S_\sigma^{-1}x_k$, where
$\widetilde x_{k,n}=\sigma^{-n/2}x_{k,n}$.
Then
\begin{equation}
\widetilde x_{k+1}
=
[\sigma I_N-(1-\sigma)J]\widetilde x_k
+(1-\sigma)\mathbf{1} u_{k+1},
\label{eq:meixner-structured-scaling}
\end{equation}
where $\mathbf{1}\in\R^N$ is the vector of ones.
Each state update can be evaluated
using $O(N)$ arithmetic operations.
Moreover,
$\kappa_2(S_\sigma)=\sigma^{-(N-1)/2}$.
\end{proposition}

\begin{proof}
By \eqref{eq:meixner-matrices},
$S_\sigma^{-1}A^{\mathrm M}S_\sigma
=\sigma I_N-(1-\sigma)J$ and
$S_\sigma^{-1}B^{\mathrm M}=(1-\sigma)\mathbf{1}$,
which gives \eqref{eq:meixner-structured-scaling}.
Since
$[J\widetilde x_k]_n
=\sum_{m=0}^{n-1}\widetilde x_{k,m}$,
define $a_0=0$ and
$a_{n+1}=a_n+\widetilde x_{k,n}$.
Then
\[
\widetilde x_{k+1,n}
=
\sigma\widetilde x_{k,n}
-(1-\sigma)a_n
+(1-\sigma)u_{k+1}.
\]
The prefix sums require a constant number of
operations per component, giving an $O(N)$ update.
Finally, the largest and smallest diagonal entries
of $S_\sigma$ are $1$ and $\sigma^{(N-1)/2}$,
respectively, so
$\kappa_2(S_\sigma)=\sigma^{-(N-1)/2}$.
\end{proof}

Standard triangular matrix-vector
multiplication with $A^{\mathrm M}$
uses $N(N+1)/2$ scalar multiplications
for the state-dependent terms,
whereas the diagonally scaled
prefix-sum update has $O(N)$
arithmetic complexity.
This comparison concerns two
implementations of the same recurrence
and does not imply an advantage over
existing Meixner-like filter realizations.
The linear arithmetic complexity
does not imply stability in finite
precision.
Indeed,
$\kappa_2(S_\sigma)=\sigma^{-(N-1)/2}$,
and no uniform floating-point error
bound is established.

\begin{example}\label{example}
In this and the following example,
we restrict the Meixner history model
to finite causal histories.
Identify the histories supported
on $\{0,\ldots,k\}$ with
$\cH_k=\C^{k+1}$, where
$h_k=\operatorname{col}(u_k,\ldots,u_0)$.
The restriction of the full Meixner
analysis to these histories is represented
by $K_k\in\C^{N\times(k+1)}$, where
$[K_k]_{n,\ell}=(1-\sigma)\sigma^\ell
m_n(\ell)$ for $0\le n<N$
and $0\le\ell\le k$.
Since $m_n$ has degree $n$ with
a nonzero leading coefficient,
$m_0,\ldots,m_{N-1}$ form a basis
of the polynomials of degree less
than $N$. Evaluation at the distinct
nodes $0,\ldots,k$ has rank
$\min\{N,k+1\}$, and the nonzero
column weights preserve this rank.
Hence $\rank K_k=\min\{N,k+1\}$.

For these histories, the transport
and insertion maps are
$S_i h=\operatorname{col}(0,h)$
and $E_{i+1}u=\operatorname{col}(u,0_{i+1})$.
The exact recurrence in
Theorem~\ref{thm:meixner-realization}
gives $K_{i+1}S_i=A^{\mathrm M}K_i$.
Thus $\row(K_{i+1}S_i)\subseteq\row K_i$,
and the backward recursion
\eqref{eq:time-varying-invariant-envelope}
yields $\cV_i=\row K_i$ on every
fixed finite horizon.
By Theorem~\ref{thm:minimal-window-augmentation},
the minimum exact realization dimension
on the full finite causal history space
at time $i$ is
$\dim\cV_i=\rank K_i$.
Thus no additional independent
state directions are required.

When the future task operators
$\cC_i$ are represented by $K_i$
for $k\le i\le t$, iterating
\eqref{eq:operator-closure} gives
$K_iS(i,k)=(A^{\mathrm M})^{i-k}K_k$
and hence
$\ker\mathscr O_{t,k}=\ker K_k$.
Taking $\mathscr R_0=\C$ and allowing
arbitrary input samples gives
$\mathscr R_k=\C^{k+1}$.
By Theorem~\ref{thm:finite-horizon-quotient-synthesis},
\[
r_k=\rank K_k=\min\{N,k+1\}.
\]
Consequently, for $k\ge N-1$,
the $N$-dimensional Meixner realization
attains the minimum state dimension
for preserving all its coefficients
over the prescribed finite horizon.
For $k<N-1$, the minimum exact
realization dimension is $k+1$,
although the prescribed analysis
has $N$ coordinates.
\end{example}

\begin{example}\label{ex:meixner-selected}
Let $K_i$ be the finite causal
Meixner analysis matrices defined
in Example~\ref{example}.
Assume $N\ge2$ and $t\ge k\ge N-1$.
For a fixed $0\le q<N$, prescribe
only the Meixner coefficient of order $q$,
$K_i^{(q)}=e_{q+1}^{\T}K_i$,
$k\le i\le t$, where $e_{q+1}$
is the $(q+1)$th standard basis vector.
Denote the corresponding backward
invariant envelope by $\cV_i^{(q)}$.
The envelope recursion and
$K_{i+1}S_i=A^{\mathrm M}K_i$ give
\[
\cV_k^{(q)}
=
\row\operatorname{col}\bigl(
e_{q+1}^{\T}K_k,\,
e_{q+1}^{\T}A^{\mathrm M}K_k,\ldots,
e_{q+1}^{\T}(A^{\mathrm M})^{t-k}K_k
\bigr).
\]

By \eqref{eq:meixner-matrices},
$S_\sigma^{-1}A^{\mathrm M}S_\sigma
=\sigma I_N-(1-\sigma)J$.
For $0\le j\le q$, the rightmost
nonzero entry of $e_{q+1}^{\T}J^j$
occurs in column $q+1-j$.
These rows are linearly independent,
whereas $e_{q+1}^{\T}J^j=0$
for $j>q$.
The binomial expansion and invertible
diagonal scaling therefore imply
\[
\rank\operatorname{col}\bigl(
e_{q+1}^{\T},
e_{q+1}^{\T}A^{\mathrm M},\ldots,
e_{q+1}^{\T}(A^{\mathrm M})^{t-k}
\bigr)
=\min\{q+1,t-k+1\}.
\]
Since $\rank K_k=N$, right
multiplication by $K_k$ preserves
the rank of matrices with $N$
columns. Moreover,
$\rank K_k^{(q)}=1$.
Consequently,
\[
\dim\cV_k^{(q)}
=\min\{q+1,t-k+1\},
\qquad
\dim\cV_k^{(q)}
-\rank K_k^{(q)}
=\min\{q,t-k\}.
\]
By Theorem~\ref{thm:minimal-window-augmentation},
the latter quantity is the minimum
number of additional independent
state directions. For $k<t$,
it also equals the restricted
propagated closure defect rank in
Corollary~\ref{cor:propagated-closure-defect-rank}.
Thus the minimum augmentation depends
on both the prescribed coefficient
and the remaining horizon.
The zeroth coefficient requires none,
whereas the highest coefficient
requires $N-1$ additional directions
when $t-k\ge N-1$.
In contrast, the full Meixner
analysis requires no augmentation.
\end{example}

\section{Numerical studies}
\label{sec:numerical}

We examine exact recurrences, noncommuting
transports, Gram conditions, and minimum
realizations. The results concern the stated
numerical settings, rather than general
performance guarantees. Unless stated
otherwise, calculations use IEEE 754
double precision.

\subsection{Structured realizations and Gram bounds}
\label{subsec:structured-numerics}

For the Meixner realization in
Theorem~\ref{thm:meixner-realization},
take $\sigma=0.92$ and $512$ independent
standard Gaussian inputs generated by
NumPy PCG64 with seed $7$, after
discarding $256$ draws.
At $N=6$, the maximum infinity norm
difference between direct evaluation
of \eqref{eq:meixner-state} and the
recurrence is $7.22\times10^{-16}$.

We compare direct history summation,
the matrix-vector recurrence using
$A^{\mathrm M}$ in
\eqref{eq:meixner-matrices},
and the scaled recurrence in
\eqref{eq:meixner-structured-scaling}
evaluated using prefix sums.
The tests use $N\in\{6,32,128\}$
and the same input sequence.
For direct summation, the coefficients
$(A^{\mathrm M})^\ell B^{\mathrm M}$
are precomputed using
\eqref{eq:meixner-structured-scaling},
and their computation is excluded
from the timings.
Table~\ref{tab:meixner-update-comparison}
reports median computation times over
seven runs after one preliminary run.
The last column gives the maximum
absolute error among the two recurrences,
compared with direct summation over all
state components and time indices.
Timings were measured on an AMD EPYC 9V74
processor using NumPy 2.3.5 with one
OpenBLAS thread.

\begin{table}[htbp]
\footnotesize
\caption{Meixner update times in milliseconds
and maximum absolute errors in the stated
environment.}
\label{tab:meixner-update-comparison}
\begin{center}
\setlength{\tabcolsep}{4pt}
\begin{tabular}{|c|c|c|c|c|} \hline
$N$ & \textbf{Direct sum} & \textbf{Matrix}
& \textbf{Prefix sum} & \textbf{Max. error} \\ \hline
$6$ & $1.45$ & $1.21$ & $2.52$ & $9.99\times10^{-16}$ \\
$32$ & $4.37$ & $1.21$ & $2.68$ & $1.28\times10^{-15}$ \\
$128$ & $16.46$ & $1.95$ & $3.10$ & $1.28\times10^{-15}$ \\ \hline
\end{tabular}
\end{center}
\end{table}

Although the prefix sum update has
$O(N)$ arithmetic cost, it is not
faster in this implementation.
Since the reference coefficients
are generated by the structured
recurrence, these errors compare
implementations but do not independently
establish accuracy at larger $N$.
The separate evaluation of
\eqref{eq:meixner-state} at $N=6$
provides an independent recurrence check.
No advantage over existing Meixner
filter realizations or uniform
floating point error bound is established.

We next test the noncommuting matrix mode
realization of Theorem~\ref{thm:matrix-mode-dho}.
Take $N=4$, $\sigma=0.94$, and
$b_k=b=[1,-0.4,0.7,0.25]^{\T}$.
For all $k\in\mathbb Z$, define
\[
R_k=U_kD(V_k^{\mathrm{rot}})^{\T},\qquad
D=\operatorname{diag}(0.76,0.84,0.91,0.96),
\]
where
\[
\begin{aligned}
U_k&=Q_*\operatorname{diag}\bigl(
\operatorname{Rot}(0.31\sin(0.071k)),
\operatorname{Rot}(0.27\cos(0.053k))\bigr),\\
V_k^{\mathrm{rot}}&=Q_*\operatorname{diag}\bigl(
\operatorname{Rot}(0.23\cos(0.043k+0.3)),
\operatorname{Rot}(-0.29\sin(0.061k))\bigr).
\end{aligned}
\]
Here
$\operatorname{Rot}(\theta)=
\begin{bmatrix}
\cos\theta&-\sin\theta\\
\sin\theta&\cos\theta
\end{bmatrix}$
and
\[
Q_*=\begin{bmatrix}
-a&0&0.22a&0\\
0&-d_*&0&-0.18d_*\\
0.22a&0&a&0\\
0&-0.18d_*&0&d_*
\end{bmatrix},
\]
with $a=(1+0.22^2)^{-1/2}$ and
$d_*=(1+0.18^2)^{-1/2}$.
Use $u_0,\ldots,u_{600}$ from NumPy PCG64
seed $811$, with independent standard
Gaussian entries and $u_j=0$ for $j<0$.
Over the tested indices, the norms
$\norm{R_{k+1}R_k-R_kR_{k+1}}$
range from $2.79\times10^{-4}$ to
$2.08\times10^{-3}$. The maximum infinity
norm difference between direct analysis
and \eqref{eq:matrix-mode-matrices} is
$6.3\times10^{-15}$.

For $k=0,\ldots,600$, Gram sums with
$\ell=0,\ldots,700$ have minimum and maximum
eigenvalues $2.18\times10^{-5}$ and $0.3614$.
Since $\norm{R_k}=0.96$, their theoretical
series truncation error is below
$1.56\times10^{-44}$, excluding rounding
errors. We examine the state transition
bounds in
\eqref{eq:matrix-mode-transition-bounds}.
At $k=600$, the decimal logarithms
of the largest and smallest singular values
of $\cF(k,0)$ are $-32.56$ and $-77.28$,
compared with the upper and lower bounds
$-31.39$ and $-78.96$, respectively.
The bounds hold over the tested indices
using the corresponding Gram eigenvalue
extrema. The extreme singular values are
evaluated using separately rescaled
ordered products and inverse products.
Calculations with $100$ decimal digits,
using the same computed Gram normalization
matrices, agree at $k=100,300,600$.

\begin{figure}[htbp]
\centering
\includegraphics[width=0.75\textwidth]{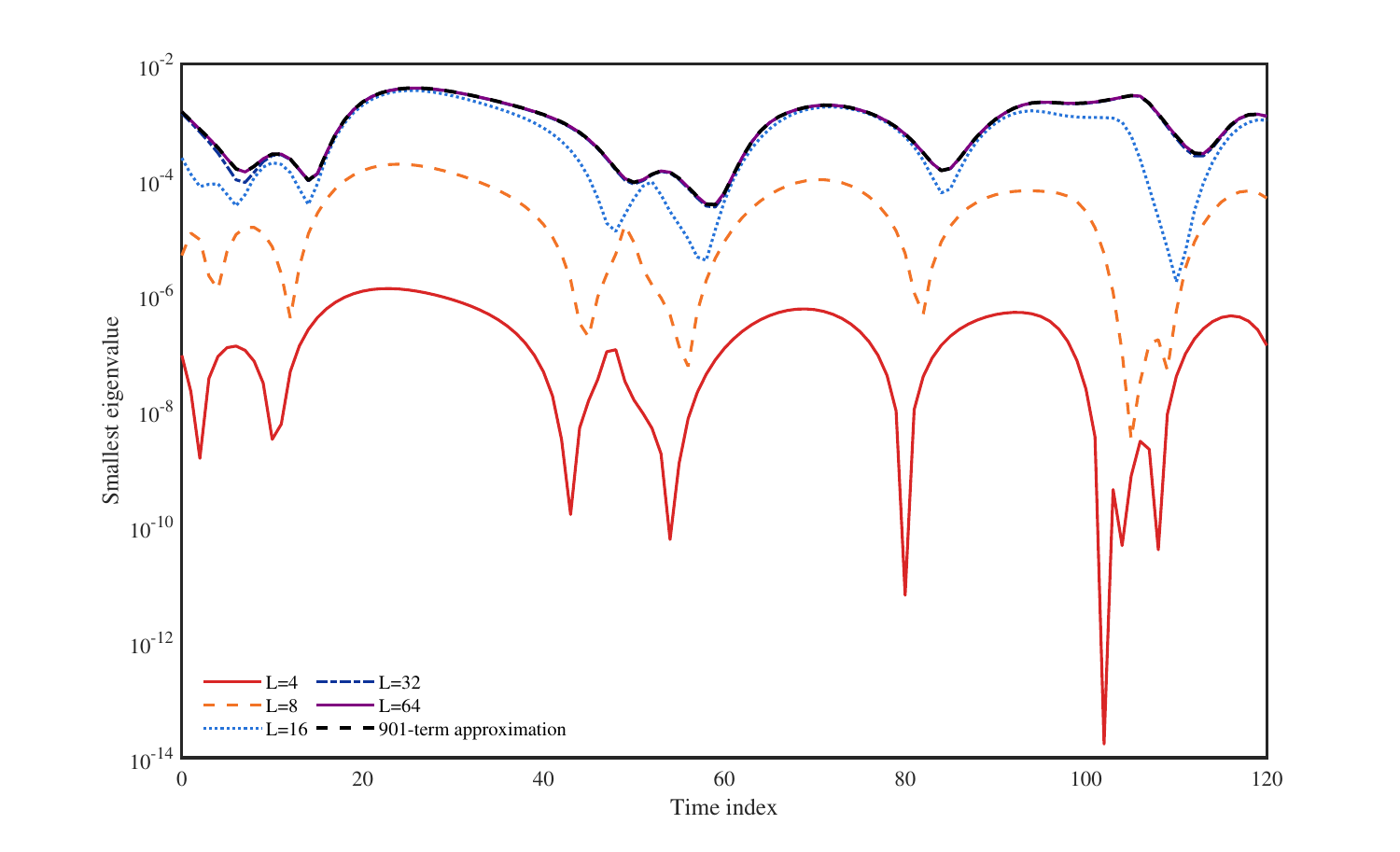}
\caption{Smallest eigenvalues of finite lag
Gram matrices and approximations of the
infinite Gram matrices using $901$ terms
for $k=0,\ldots,120$.}
\label{fig:finite-gram-certificate}
\end{figure}

We next examine the finite Gram condition
of Theorem~\ref{thm:finite-gram-certificate}.
Figure~\ref{fig:finite-gram-certificate}
compares the smallest eigenvalues of
$W_{k,L}$ for $L\in\{4,8,16,32,64\}$
and Gram sums with $\ell=0,\ldots,900$
at $k=0,\ldots,120$. The minimum eigenvalue
of $W_{k,32}$ is $3.37\times10^{-5}$,
compared with $3.71\times10^{-5}$ for
the $901$ term Gram sums.
The latter have maximum eigenvalue $0.348$,
with a theoretical truncation error
below $5.33\times10^{-57}$.
The results illustrate the finite Gram
lower bound over the tested indices.

We finally examine the common metric condition
of Proposition~\ref{prop:metric-gram-certificate}
for a period two model with $N=3$,
$\sigma=0.9$, and $L=12$. Define
\[
T=\begin{bmatrix}
1&3&-0.4\\0&1&2\\0&0&1
\end{bmatrix},\qquad
Q=\begin{bmatrix}
0.8&-0.6&0\\0.6&0.8&0\\0&0&1
\end{bmatrix},
\]
and, for every $j\in\mathbb Z$,
\[
R_{2j}=T^{-1}\operatorname{diag}(0.68,0.80,0.86)T,
\qquad
R_{2j+1}=T^{-1}Q\operatorname{diag}(0.74,0.83,0.85)Q^{\T}T.
\]
Take $b_i=b=T^{-1}[1,-0.4,0.7]^{\T}$
and $P=T^{\T}T/\operatorname{tr}(T^{\T}T)$.
The Euclidean hypothesis of
Theorem~\ref{thm:finite-gram-certificate}
fails because
$\sigma\max_{j=0,1}\norm{R_j}^2
\approx1.166009>1$.
In contrast,
$\max_{j=0,1}\norm{P^{1/2}R_jP^{-1/2}}^2
\approx0.7396<r_P^2=0.75$,
$\sigma r_P^2=0.675<1$, and
$b^{\T}Pb=165/1616<\beta_P^2=0.102104$.
Using $g_j(0)=b$ and the ordered
products
$g_j(\ell)=R_jR_{j-1}\cdots
R_{j-\ell+1}b$ for $\ell\ge1$,
exact rational arithmetic shows that,
for both $j=0,1$, the leading principal
minors of $W_{j,12}-2\times10^{-6}P^{-1}$
of orders $1,2,3$ exceed $11$, $10^{-2}$,
and $7\times10^{-9}$, respectively.
The smallest eigenvalues of
$P^{1/2}W_{j,12}P^{1/2}$ at $j=0,1$
are approximately $2.20\times10^{-6}$
and $3.79\times10^{-6}$.
Since $W_{k+2,12}=W_{k,12}$,
the finite Gram lower bound holds
for all $k\ge0$. Consequently,
Proposition~\ref{prop:metric-gram-certificate}
gives
\[
2\times10^{-6}P^{-1}
\preceq G_k
\preceq 0.03142P^{-1},
\qquad k\ge0.
\]
Thus uniform Gram bounds hold
although the Euclidean sufficient
condition is not satisfied.

\subsection{Cross-projection and minimal realizations}
\label{subsec:core-numerics}

Consider the sliding history
$h_k=[u_{k-11},\ldots,u_k]^{\T}\in\R^{12}$.
Let $K_k=K$, $S_k=S$, and $E_{k+1}=E$,
where $S_{i,i+1}=1$ for $1\le i<12$,
all other entries of $S$ vanish,
$E=e_{12}$, and
$K=\frac12[0_{1\times8},1,1,1,1]$.
Since $KK^{\T}=1$, the cross-projection
update has $A_k^{\mathrm{cp}}=KSK^{\T}=3/4$
and $B_{k+1}=KE=1/2$.
Its defect is
$\Delta_k^{\mathrm{cp}}
=KS(I_{12}-K^{\T}K)\ne0$.

For $t-k\ge3$, the rows of
$K,KS,KS^2,KS^3$ are independent,
while $KS^4=0$.
Indeed, the submatrix of
$\operatorname{col}(K,KS,KS^2,KS^3)$
formed by columns $9,\ldots,12$ has
determinant $1/16$.
Consequently, $\dim\cV_k=4$ and
\[
\rank\bigl(\mathsf D_{t,k}^{\mathrm{cl}}
\big|_{\ker K}\bigr)
=\dim\cV_k-\rank K=3.
\]
Thus, for $t-k\ge3$,
Theorem~\ref{thm:minimal-window-augmentation}
gives a minimum exact state dimension
of four within the admissible class
defined by
\eqref{eq:augmentation-requirements}.
This minimum is attained by retaining
the four most recent input samples.

\begin{figure}[htbp]
\centering
\includegraphics[width=0.75\textwidth]{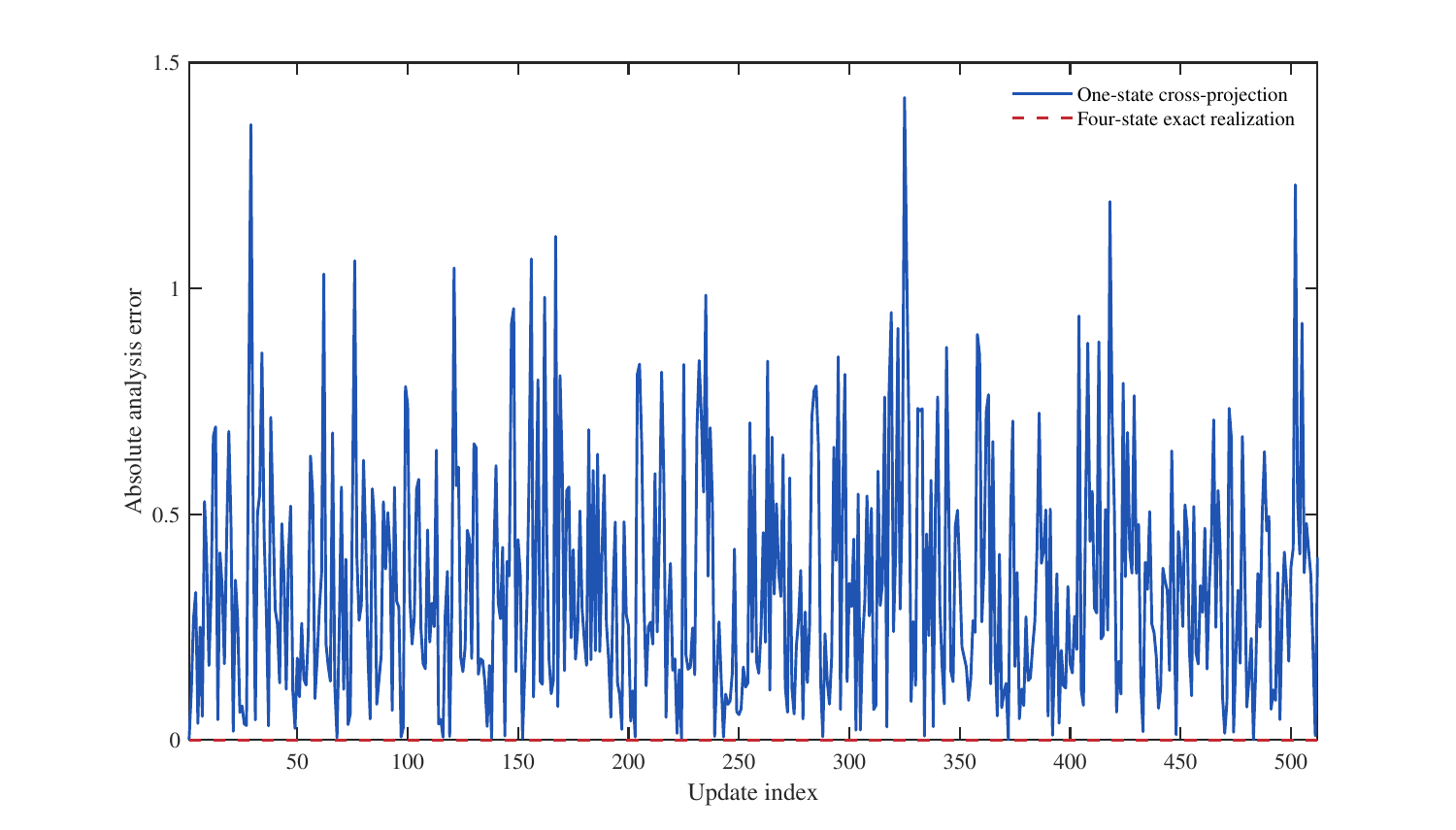}
\caption{Absolute analysis errors for one
seeded input sequence relative to direct
history analysis. The realization with
four states is exact up to rounding.}
\label{fig:minimal-augmentation-comparison}
\end{figure}

With zero initial history, use $512$
independent standard Gaussian inputs
from NumPy PCG64 seed $20260812$.
The scalar cross-projection recurrence
has root mean square analysis error
$0.434643$ and maximum absolute error
$1.422794$. The four state realization
agrees with direct history analysis to
rounding accuracy, as illustrated in
Figure~\ref{fig:minimal-augmentation-comparison}.

We next examine minimal realizations
on reachable histories through
\eqref{eq:reachable-hankel-rank},
using $50$ random systems.
In each trial, draw uniformly and
independently the horizon
$t\in\{4,\ldots,8\}$, history dimension
$d\in\{6,\ldots,10\}$, task dimension
$p\in\{1,2,3\}$, and input dimension
$m\in\{1,2,3\}$.
These dimensions remain fixed
within each trial.
For $i=0,\ldots,t-1$,
obtain $\widetilde Q_i$ by unpivoted
QR factorization of a new
$d\times d$ standard Gaussian matrix,
draw the entries of $\lambda_i$
independently from the uniform
distribution on $[0.45,0.92]$, and set
$S_i=\widetilde Q_i
\operatorname{diag}(\lambda_i)
\widetilde Q_i^{\T}+0.03Z_i$,
where $Z_i$ has independent
standard Gaussian entries.
Next draw $E_0,\ldots,E_t$ of size
$d\times m$, followed by
$\cC_0,\ldots,\cC_t$ of size
$p\times d$, with independent
standard Gaussian entries.

Set $\mathscr R_0=\operatorname{range}E_0$.
Following the reachable history recursion
\eqref{eq:reachable-history-recursion},
construct the reachable maps from
$\mathscr G_{0,0}=E_0$ and
$\mathscr G_{k+1,0}
=[S_k\mathscr G_{k,0},E_{k+1}]$.
All draws use a single NumPy PCG64
generator initialized with seed
$20260801$, in the stated order.

Numerical ranks count singular values
exceeding $10^{-10}$ times the largest,
with zero matrices assigned rank zero.
At all $360$ tested cuts, the computed
quotient dimensions agree with the
numerical ranks of the corresponding
reachable future-task matrix, consistent
with \eqref{eq:reachable-hankel-rank}.
Across these reachable future-task
matrices, the
smallest retained singular value ratio
is $4.70\times10^{-4}$, while the
largest discarded ratio is
$1.83\times10^{-16}$.

\section{Conclusions}
\label{sec:conclusions}

This paper establishes exact closure
conditions for prescribed discrete
history analyses and distinguishes
two finite-horizon minimality problems.
Reachable task quotients determine
minimum state dimensions for preserving
specified future tasks on reachable
histories, whereas
backward invariant envelopes characterize
minimum exact realizations on full
history spaces within the admissible class.
The minimum number of additional
independent state directions equals
the rank of the propagated closure
defect operator restricted to the
current analysis kernel.
For countable histories, matrix mode
transport yields structured exact
recurrences, uniform exponential decay
estimates, and sufficient conditions
for uniform Gram bounds.
The Meixner specialization provides
explicit recurrence matrices, an
$O(N)$ arithmetic state update,
and finite-horizon minimum augmentation
dimensions for selected coefficients.
Numerical experiments illustrate
these results.
Future work includes stochastic
extensions and further analysis of
numerical conditioning.

\end{document}